\documentclass[12pt]{article}       

\usepackage{amsmath,amsthm}
\usepackage{amssymb}
\usepackage{latexsym}
\begin{document}

\newcommand{\comment}[1]{}    
\newcommand{\hs}{\enspace}
\newcommand{\hhs}{\thinspace}
\newcommand{\real}{\ifmmode {\rm R} \else ${\rm R}$ \fi}

\newtheorem{theorem}{Theorem}
\newtheorem{lemma}[theorem]{Lemma}         
\newtheorem{corollary}[theorem]{Corollary}
\newtheorem{definition}[theorem]{Definition}
\newtheorem{claim}[theorem]{Claim}
\newtheorem{conjecture}[theorem]{Conjecture}
\newtheorem{proposition}[theorem]{Proposition}

\newtheorem{remark}[theorem]{Remark}



\def\edge{\leftrightarrow}
\def\noedge{\not\leftrightarrow}
\def\twoedge{\Leftrightarrow}
\def\to{\rightarrow}
\def\Hrl{H^{(r)}_{l+1}}
\def\Krl{{\cal K}^{(r)}_l}
\def\Krl{{\cal K}^{(r)}_{l+1}}
\def\cF{{\cal F}}
\def\cG{{\cal G}}
\def\cH{{\cal H}}
\def\e{\varepsilon}
\def\bF{
{\cal {\bf F}}}
\def\odel{o_{\delta}}
\def\od1{o_{\delta}(1)}
\def\bF{{\bf F}}
\def\oe1{o_{\varepsilon}(1)}


\title{\bf Counting substructures II:  triple systems}
\author{Dhruv Mubayi
\thanks{ Department of Mathematics, Statistics, and Computer
Science, University of Illinois, Chicago, IL 60607.  email:
mubayi@math.uic.edu; research  supported in part by  NSF grant DMS
0653946
}}
\date{\today}
\maketitle

\begin{abstract}
For various triple systems  $F$, we give tight lower
bounds on the number of copies of $F$ in  a triple system with a
prescribed number of vertices and edges.  These are the first such
results for hypergraphs, and extend earlier theorems of Bollob\'as,
Frankl, F\"uredi, Keevash, Pikhurko, Simonovits, and Sudakov who proved
that there is one copy of $F$.

A sample result is the following:  F\"uredi-Simonovits \cite{FS} and independently
Keevash-Sudakov \cite{KSFano} settled an old conjecture of S\'os \cite{Sos} by
proving that the maximum number of
triples in an $n$ vertex triple system (for $n$ sufficiently large)
that contains no copy of the Fano plane is
$p(n)={\lceil n/2\rceil \choose 2}\lfloor n/2 \rfloor+{\lfloor
n/2\rfloor \choose 2}\lceil n/2 \rceil.$

We prove that there is an absolute constant $c$ such that if $n$ is sufficiently large and $1 \le q \le cn^2$, then every $n$ vertex  triple
system with $p(n)+q$ edges contains at least
$$6q\left({\lfloor n/2 \rfloor \choose 4}+(\lceil n/2 \rceil-3){\lfloor n/2 \rfloor\choose 3}\right)$$
copies of the Fano plane. This is sharp for $q\le n/2-2$.


Our proofs use  the recently proved hypergraph removal lemma and
stability results for the corresponding Tur\'an problem.

\end{abstract}

\section{Introduction}
Many mathematical problems enjoy the supersaturation
phenomenon which, broadly speaking, says that once we have
sufficiently many objects of a particular type to guarantee the
existence of a specific configuration, then we find not just one
but many copies of this configuration.  The objects can
be edges in a graph, points in the plane, subsets of integers, etc.
One well-known example is Szemer\'edi's theorem about the existence
of arithmetic progressions of length $k$ in a subset $S \subset \{1,
\ldots, n\}$ with $S$ sufficiently large. It is known that if $|S|\ge \e
n$ (with $\e>0$ fixed and $n$ sufficiently large) then we are
guaranteed not just one, but $c_{k, \e}n^2$ arithmetic progressions
of length $k$ from $S$ (see Tao \cite{T2} or Varnavides \cite{V}).

Perhaps the early examples of this phenomenon came from graph
theory. Mantel proved that  a graph with $n$ vertices and
$\lfloor n^2/4 \rfloor+1$ edges contains a triangle. Rademacher
extended this by showing that there are   at least $\lfloor n/2
\rfloor$ copies of a triangle.  Subsequently, Erd\H os \cite{E1, E2}
proved that if $q <cn$ for some small constant $c$, then $\lfloor
n^2/4 \rfloor+q$ edges guarantees at least $q\lfloor n/2 \rfloor$
triangles.  Later Lov\'asz and Simonovits \cite{LS} proved that the
same statement holds with $c=1/2$, thus confirming an old conjecture
of Erd\H os.  They also proved similar results for complete graphs.

In this paper (the second in a  series) we initiate the study of this phenomenon to $k$-uniform hypergraphs ($k$-graphs for short).
In the first paper of this series \cite{KaMu}, we had extended the results of Erd\H os and Lov\'asz-Simonovits in two ways. First,
 we proved such statements for the broader class of color critical graphs.  Second,
 we showed that all the copies of the required subgraph were incident to a small number of edges or vertices. For example, in a graph with $n$ vertices and $\lfloor n^2/4 \rfloor+q$ edges, \cite{E1, E2, LS} do not give information about how the $q\lfloor n/2 \rfloor$ triangles are distributed. In \cite{KaMu}, we proved that as long as $q=o(n)$ there are  $(1-o(1))qn/2$ triangles incident with at most $q$ vertices.

The main new tool we have at our disposal is the recently proved hypergraph removal lemma, which is a
consequence of the hypergraph regularity lemma
 (see Gowers \cite{G}, Nagle-R\"odl-Schacht \cite{NRS}, R\"odl-Skokan \cite{RS}, Tao \cite{T}).
 The novelty in this project is the use of the removal lemma to count
 substructures in hypergraphs rather precisely.

\begin{theorem} {\bf (Hypergraph Removal Lemma \cite{G, NRS, RS, T})} \label{removal}
Fix $k \ge 2$ and a $k$-graph $F$ with $f$ vertices. Suppose that an
$n$ vertex $k$-graph $\cH$ has at most $o(n^f)$ copies of $F$. Then
there is a set of edges in $\cH$ of size $o(n^k)$ whose removal from
$\cH$ results in a $k$-graph with no copies of $F$.
\end{theorem}

Given a $k$-graph $F$, let ex$(n, F)$, the Tur\'an number of $F$, be the maximum number of edges in an $n$
vertex $k$-graph with no copy of $F$. For $k>2$,  determining
the  Tur\'an number is a very difficult problem, and there are only sporadic results. Many of these
were obtained recently by using the so-called stability approach first introduced by Erd\H os and
Simonovits \cite{S} in the late 1960's. Here we take this project one step further by giving asymptotically
 sharp results on the number of copies of a $k$-graph $F$ in a $k$-graph with $n$  vertices and ex$(n,F)+q$
 edges. In two cases we are able to count the exact minimum number of copies even though this number is quite complicated (see the abstract).

 In essentially all cases where ex$(n,F)$ is known (when $k=3$), it turns out that one is guaranteed many copies
 of $F$ as long as there are ex$(n,F)+1$ edges, so we extend all
 previous results that determine ex$(n,F)$.
It is somewhat surprising that although determining ex$(n,F)$ for
these hypergraphs $F$ is quite difficult (in some cases they were
decades old conjectures that were only recently settled), we are
able to count quite precisely the number of copies of $F$ as long as
the  number of extra edges $q$ is not too large. Typically we can
allow $q=o(n^2)$ for the 3-graphs we consider.

Here we should also mention the relationship between this project
and recent work of Nikiforov \cite{N} and Razborov \cite{R} that
gives asymptotically sharp estimates on the minimum number of
triangles in a graph with $n$ vertices and $\lfloor n^2/4 \rfloor
+q$ edges, where $q=\Omega(n^2)$.  There are at present no such
results for $k$-graphs for $k>2$, and little hope of achieving them.
Moreover, even if such results were to be proved, they would apply
only  when $q=\Omega(n^k)$, so the results of the type \cite{N, R}
will not overlap with ours.

Our proofs all have the following basic structure: Suppose we are given $\cH$ with sufficiently many edges and we
wish to find many copies of $F$ in $\cH$. First we observe that if the number of copies of $F$ is very large, then
we already have the bound sought.  Consequently, we can use the hypergraph removal lemma to delete a small proportion
of edges of $\cH$ so that the resulting triple system has no copies of $F$.  Next we use the stability results that
guarantee the approximate structure of $\cH$.  At this point the techniques depend highly on the particular
structure of $F$ and of $\cH$. The technical details are  more involved than for the usual Tur\'an problem,
 since it is not enough to find just one copy of $F$.
 At the end of the analysis, we are able to describe quite precisely how the copies of $F$ are distributed
 within $\cH$.

We illustrate our approach on four excluded hypergraph problems, each of which has been studied quite a lot.

\begin{definition} Let  $F$ have the property that for sufficiently large $n$, there is a unique (up to isomorphism) 3-graph $\cH(n)$ with ex$(n,F)$ edges. Let $c(n,F)$ be the minimum number of copies of $F$ in the 3-graph obtained from $\cH(n)$ by adding an edge, where the minimum is taken over all possible ways to add an edge.
\end{definition}

Our theorems all say that if $\cH$ is an $n$ vertex 3-graph with ex$(n,F)+q$ edges, then the number of copies of $F$ in $\cH$ is essentially at least $qc(n,F)$.
In the next subsections we will  state our results precisely.

There remain a few more exact 4-graph results in the literature where we could possibly use this approach for
the counting problem. We will give the corresponding counting results for all of these problems in a forthcoming paper \cite{M3}, the third in this series.

Notation: We associate a hypergraph with its edge set. The number of edges in a hypergraph $\cH$ is $|\cH|$. Given hypergraphs
$F, \cH$ ($F$ has $f$ vertices), a copy of $F$ in $\cH$ is a subset of
$f$ vertices and $|F|$ edges of $\cH$ such that the subhypergraph formed
by this set of vertices and edges is isomorphic to $F$.  In other words, if we denote $Aut(F)$ to be the number of automorphisms of $F$, then the number of copies of $F$ in $\cH$ is the number of edge-preserving injections from $V(F)$ to $V(\cH)$ divided by $Aut(F)$.
For a set $S$ of vertices, define $d_{\cal H}(S)$ to be the number of edges of $\cH$ containing $S$.  If $S=\{v\}$, we simply write $d_{\cal H}(v)$.
We will omit floor and ceiling symbols whenever they are not crucial, so that the presentation is clearer.

\subsection{Fano plane}

Let ${\bf F}$ be the projective plane of order two over the finite field of order two.
An explicit description of ${\bf F}$ is $\{124, 235, 346, 457, 561, 672, 713\}$, obtained from the difference set $\{1,2,4\}$ over $Z_7$. It is well known that ${\bf F}$ is not 2-colorable, hence it cannot be a subgraph of any 2-colorable 3-graph. Say that
a 3-graph $\cH$ is bipartite (or 2-colorable) if it has
a vertex partition $A \cup B$ such that every edge intersects both
parts.  Let $P^3(n)$ be the bipartite 3-graph with the maximum
number of edges. Note that
$$p^3(n):=|P^3(n)|=\max_a\left\{ {a \choose 2}(n-a)+{n-a \choose
2}a\right\} =(3/4+o(1)){n \choose 3}$$ is uniquely achieved by choosing
$a\in\{\lfloor n/2 \rfloor, \lceil n/2\rceil\}$.

 S\'os \cite{Sos} conjectured, and Keevash-Sudakov \cite{KSFano} and F\"uredi-Simonovits \cite{FS} independently proved that among all $n$ vertex 3-graphs ($n$ sufficiently large) containing no copy of  ${\bf F}$, the unique one with the maximum number of edges is $P^3(n)$.  Thus $c(n, \bF)$ is defined and in fact
 $$c(n,{\bf F}):=6\left({\lfloor n/2 \rfloor \choose 4}+(\lceil n/2 \rceil-3){\lfloor n/2 \rfloor\choose 3}\right)=(20+o(1))(n/4)^4.$$
 This is achieved by adding an  edge to the part of size $\lceil n/2 \rceil$.  Indeed, if we add a triple $123$ to this part, then one way to make a copy of ${\bf F}$ is to take four points $a,b,c,d$ from the other part, partition the six pairs among $\{a,b,c,d\}$  into three perfect matchings $m_1, m_2, m_3$, and use the edges $\{i\} \cup p$ where $p \in m_i$, for each $i$ to form a copy of ${\bf F}$.  There are  ${\lfloor n/2\rfloor \choose 4}$ ways to pick $a,b,c,d$ and for each such choice there are six ways to choose $m_1, m_2, m_3$. The only other way to form a copy of ${\bf F}$ using $123$ is to pick four points $a,b,c,d$ with $a$ in the same part as $1$ and $b,c,d$ in the other part. Then proceeding as before, we obtain $6(\lceil n/2\rceil-3){\lfloor n/2\rfloor \choose 3}$ copies of ${\bf F}$.  Altogether we obtain $c(n, {\bf F})$ copies.
 
 Our first result shows that a 3-graph with $p^3(n)+q$ edges has at least as many copies
  of $\bF$ as a 3-graph obtained from $P^3(n)$ by adding $q$ edges in an optimal way.  The precise number we can add is
 $$q(n, \bF)=\begin{cases}
 n & \hbox{ if $n$ is even and $n/2 \equiv 0$ (mod 4)} \\
 n-2 & \hbox{ if $n$ is even and $n/2 \equiv 1$ (mod 4)} \\
 n-4 & \hbox{ if $n$ is even and $n/2 \equiv 2,3$ (mod 4)} \\
 \lceil n/2 \rceil & \hbox{ if $n$ is odd and $\lceil n/2 \rceil \equiv 0$ (mod 4)} \\
 \lceil n/2 \rceil-1 & \hbox{ if $n$ is odd and $\lceil n/2 \rceil \equiv 1$ (mod 4)} \\
 \lceil n/2 \rceil-2 & \hbox{ if $n$ is odd and $\lceil n/2 \rceil \equiv 2,3$ (mod 4).} \\ \end{cases}$$
 \medskip

\begin{theorem} \label {Fanoq}
There exists an absolute constant $\e>0$ such that if $n$ is sufficiently large and  $1 \le q \le \e n^2$, then the following holds:

$\bullet$ Every $n$ vertex 3-graph with $p^3(n)+q$ edges contains at least $qc(n,\bF)$ copies of $\bF$.  This is sharp for all $q \le q(n ,\bF)$.

$\bullet$ If $q>q(n,\bF)$, then every $n$ vertex 3-graph with $p^3(n)+q$ edges contains at least $qc(n,\bF)+1$ copies of $\bF$.
\end{theorem}

{\bf Remark.} For $q>q(n, \bF)$, our proof actually gives at least $qc(n, \bF)+2{\lfloor n/2 \rfloor \choose 2}$ copies of $\bF$.

To see that Theorem \ref{Fanoq} is tight for $q\le q(n ,\bF)$ observe that we may add $q$  edges to $P^3(n)$ with every two edges sharing zero or two points. If $n$ is even, we do this by adding to each part of $P^3(n)$ the maximum number of edge that pairwise share zero or two points. This is achieved by adding disjoint copies of $K_4^3$, the complete 3-graph on four points, or collections of edges that pairwise share the same two points.
If $n$ is odd, then we add edges only to the larger part. Each added edge lies in exactly $c(n,{\bf F})$ copies of ${\bf F}$ and no two added edges lie in any copy of $\bF$. So the total number of copies of $\bF$ is exactly $qc(n, \bF)$.

Theorem \ref{Fanoq} is asymptotically sharp in a much larger range of $q$. In particular, we have the following.

\begin{proposition} \label{pasch}
For every $\e>0$ there exists $\delta>0$ and $n_0$ such that the
following holds for all $n>n_0$ and $q<\delta n^2$. There is an $n$ vertex 3-graph with $p^3(n)+q$
edges and at most $(1+\e)qc(n,\bF)$ copies of $\bF$.
\end{proposition}

\subsection{Cancellative triple-systems}

Say that a 3-graph is cancellative if whenever $A \cup B=A \cup C$ we have $B=C$. An equivalent definition
is to simply say that the 3-graph does not contain a copy of two particular 3-graphs:
$F_5=\{123, 124, 345\}$ and $K_4^{3-}=\{123, 124, 234\}$.
Write
$$t^3(n)=\left\lfloor {n\over 3} \right\rfloor \left\lfloor {n+1\over 3} \right\rfloor \left\lfloor {n+2\over 3}
\right \rfloor$$ for the number of edges in $T^3(n)$, the complete
3-partite 3-graph with the maximum number of edges.  It is easy to see that $T^3(n)$ is cancellative.

 Katona conjectured,  and Bollob\'as \cite{B} proved, that the maximum number of edges in an $n$ vertex cancellative 3-graph
is $t^3(n)$, and equality holds only for $T^3(n)$. Later Frankl and F\"uredi \cite{FF} refined this by proving
the same result (for $n>3000$) even if we just forbid $F_5$.  Recently, Keevash and the author \cite{KM}
gave a new proof of the Frankl-F\"uredi result while reducing the smallest $n$ value to 33.

It is easy to see that $c(n, F_5)=3(n/3)^2+\Theta(n)$ and this is achieved by adding
a triple to $T^3(n)$  with two points in the largest part. In fact, even if we add a triple within one of the parts  we get almost the same number of copies of $F_5$.  Our second result shows that this is optimal, even when we are allowed to add as many as $o(n)$ edges.

\begin{theorem} \label {F5q}
For every $\e>0$ there exists $\delta>0$ and $n_0$ such that the
following holds for $n>n_0$. Let $\cH$ be a 3-graph with $t^3(n)+q$
edges where $q<\delta n$.  Then the number of copies of $F_5$ in
$\cH$ is at least $q(1-\e)c(n, F_5)$. This is asymptotically sharp for
$1\le q<\delta n$. Moreover, if the number of copies is less than
$\delta n^3$, then there is a collection of $q$ distinct edges that
each lie in $(1-\e)c(n, F_5)$ copies of $F_5$ with no two of these edges accounting for the same copy of $F_5$.
\end{theorem}

\subsection{Independent neighborhoods}

The neighborhood of a $(k-1)$-set $S$ of vertices in a $k$-graph is the set of vertices $v$ whose union with $S$ forms an edge. A set is independent if it contains no edge. We can rephrase Mantel's theorem as follows: the maximum number of edges in a 2-graph with independent neighborhoods is $\lfloor n^2/4 \rfloor$. This formulation can be generalized to $k>2$ and there has been quite a lot of recent activity on this question. We focus here on $k=3$, and observe that a 3-graph has independent neighborhoods if and only if it contains no copy of
$B_5=\{123, 124, 125, 345\}$. A 3-graph $\cH$ has a
$(2,1)$-partition if it has a vertex partition $A \cup B$ such that
$|e \cap A|=2$ for all $e \in \cH$. Let $B^3(n)$ be the 3-graph with
the maximum number of edges among all those that have $n$ vertices
and a $(2,1)$-partition. Note that
$$b^3(n):=|B^3(n)|=\max_a {a \choose 2}(n-a)=(4/9+o(1)){n \choose 3}$$
is achieved by choosing $a=\lfloor 2n/3 \rfloor$ or $a=\lceil 2n/3
\rceil$.

The author and R\"odl \cite{MR} conjectured, and F\"uredi, Pikhurko, and Simonovits \cite{FPS} proved,
that among all $n$ vertex 3-graphs ($n$ sufficiently large) containing no copy of  $B_5$, the unique one with the maximum number of edges is $B^3(n)$.

It is easy to see that $c(n, B_5)=2(n/3)^2+\Theta(n)$ and this is achieved by adding
a triple to $T^3(n)$  contained in the larger part. In fact, even if we add a triple within the smaller part  we get almost the same number of copies of $B_5$.  Our third result shows that this is optimal, even when we are allowed to add as many as $o(n^2)$ edges.

\begin{theorem} \label {B5q}
For every $\e>0$ there exists $\delta>0$ and $n_0$ such that the
following holds for $n>n_0$. Let $\cH$ be a 3-graph with $b^3(n)+q$
edges where $q<\delta n^2$.  Then the number of copies of $B_5$ in
$\cH$ is at least $q(1-\e)c(n, B_5)$. This is asymptotically sharp for
$1\le q<\delta n^2$. Moreover, if the number of copies is less than
$\delta n^4$, then there is a collection of $q$ distinct edges that
each lie in $(1-\e)c(n,B_5)$ copies of $B_5$ with no two of these edges accounting for the same copy of $B_5$.
\end{theorem}

\subsection{Expanded Cliques}

 Let $L_r$ be the 3-graph
obtained from the complete graph $K_{r}$ by enlarging each edge
with a new vertex.  These new vertices are distinct for each edge, so $L_r$ has $l_r=r+{r \choose 2}={r+1 \choose 2}$ vertices and ${r \choose 2}$ edges. Write $T^3_r(n)$ for the complete $r$-partite 3-graph with the maximum number of edges. So $T^3_r(n)$ has vertex partition $V_1 \cup \cdots  \cup V_r$, where $n_i:=|V_i|=\lfloor (n+i-1)/r\rfloor$, and all triples with at most one point in each $V_i$.  Define
$$t^3_r(n):= |T^3_r(n)| =\sum_{S \in {[r] \choose 3}} \prod_{i \in S} n_i.$$
Every set of $r+1$ vertices in $T^3_r(n)$ contains two vertices in the same part, and these two vertices lie in no edge. Consequently, $L_{r+1} \not\subset T^3_r(n)$.

The author \cite{M} conjectured, and Pikhurko \cite{P} proved,
that among all $n$ vertex 3-graphs containing no copy of  $L_{r+1}$ ($r \ge 3$ fixed, $n$ sufficiently large), the unique one with the maximum number of edges is $T^3_r(n)$.  Thus $c(n, L_{r+1})$ is defined and in fact
$$c(n, L_{r+1})=(1+o(1))\left(\left(1-{2\over r}\right)n \right)^{{r+1 \choose 2}-1} \times\left({n \over r}\right)^{r-1}=\Theta(n^{l_{r+1}-3})$$and this is achieved by adding a triple with exactly two points in a largest part. Our  final results shows that this is asymptotically optimal, even when we are allowed to add as many as $o(n^2)$ edges.

\begin{theorem} \label {Lq} {\bf (Asymptotic Counting)}
Fix $r \ge 3$. For every $\e>0$ there exists $\delta>0$ and $n_0$ such that the
following holds for $n>n_0$. Let $\cH$ be a 3-graph with $t^3_r(n)+q$
edges where $q<\delta n^2$.  Then the number of copies of $L_{r+1}$ in
$\cH$ is at least $q(1-\e)c(n, L_{r+1})$. The expression $q$ is sharp for
$1\le q<\delta n^2$. Moreover, if the number of copies is less than
$\delta n^{{r+1 \choose 2}-1}$, then there is a collection of $q$ distinct edges that
each lie in $(1-\e)c(n, L_{r+1})$ copies of $L_{r+1}$ with no two of these edges accounting for the same copy of $L_{r+1}$.
\end{theorem}

Our next result improves the asymptotic counting result above to an exact result, with a more restricted range for $q$.

\begin{theorem} \label{exactLq} {\bf (Exact Counting)}
Fix $r \ge 3$, $q>0$ and let $n$ be sufficiently large. Every $n$ vertex  triple
system with $t^3_r(n)+q$ edges contains at least $q c(n, L_{r+1})$
copies of $L_{r+1}$.
\end{theorem}

Theorem \ref{exactLq} is clearly tight, as we may add an appropriate set of $q$ pairwise disjoint edges to  $T^3_r(n)$ such that each edge lies in exactly $c(n, L_{r+1})$ copies of $L_{r+1}$.

Throughout the paper we will frequently use the notation $\delta \ll \e$, which means that $\delta$, and any function of $\delta$ (that tends to zero with $\delta$) used in a proof is smaller than any function of $\e$ used in the proof.  It is pretty difficult to write the precise dependence between $\delta$ and $\e$ as one of the constraints comes from an application of the removal lemma.

\section{Counting Fano's}

In this section we will prove Theorem \ref{Fanoq} and Proposition \ref{pasch}.
We need some lemmas about binomial coefficients.

\begin{lemma} \label{lemma1} Let $x,y,t>0$ be integers with $x+y=n$, $t<n^2$ and $s=\lceil \sqrt{2t/(n-2)}\rceil$.  Suppose that $n$ is sufficiently large and
$${x \choose 2}y+ {y \choose 2} x \ge p^3(n)-t.$$
Then $\lfloor n/2 \rfloor-s \le  x \le \lceil n/2\rceil +s$ and if $t<(n-2)/2$, then $\lfloor n/2 \rfloor-s <  x < \lceil n/2\rceil +s$.
\end{lemma}

\proof Suppose for contradiction that $x>\lceil n/2\rceil+s$ (the upper bound on $t$ ensures that $s<2\sqrt n$ and hence $x<3n/4$). Write
$$f(x)={x \choose 2}(n-x)+{n-x \choose 2}x=\frac12(n-2)x(n-x).$$
Note that $p^3(n)=f(\lfloor n/2\rfloor)=f(\lceil n/2 \rceil)$.  Our goal therefore is to obtain the contradiction $f(x)<f(\lceil n/2 \rceil)-t$. Observe that
$$f(a+1)=f(a)-\frac12(n-2)(2a+1-n).$$
Applying this repeatedly beginning with $a=\lceil n/2 \rceil$ we obtain
\begin{align}
f(x)< f(\lceil n/2\rceil+s)&=f(\lceil n/2 \rceil)-\frac12(n-2)\sum_{a=\lceil n/2\rceil}^{\lceil n/2\rceil+s-1}(2a+1-n)\notag \\
&=f(\lceil n/2 \rceil)-\frac12 s(n-2)(s+2\lceil n/2\rceil-n).\notag \end{align}
The choice of $s$ gives
\begin{equation} \label{strict} \frac12 s(n-2)(s+2\lceil n/2\rceil-n)\ge \frac12s^2(n-2)\ge t\end{equation}
and therefore $f(x)<f(\lceil n/2 \rceil)-t$.
We conclude that $x\le \lceil n/2\rceil+s$.  Repeating this argument with $x$ replaced by $y$ gives $y\le \lceil n/2\rceil+s$ and hence $x\ge \lfloor n/2\rfloor-s$.

If $t<(n-2)/2$, and $x\ge \lceil n/2\rceil+s$, then we
 only have $f(x)\le f(\lceil n/2\rceil+s)$. However the
last inequality in (\ref{strict}) is strict (since $s\ge 1$) and we again get the same contradiction.
Therefore $x< \lceil n/2\rceil+s$ and by a similar argument, $x> \lfloor n/2\rfloor-s$.
\qed
\medskip

Recall that 
 $$c(n,{\bf F}):=6\left({\lfloor n/2 \rfloor \choose 4}+(\lceil n/2 \rceil-3){\lfloor n/2 \rfloor\choose 3}\right).$$
\smallskip

\begin{lemma} \label{lemma2} Let $x,y,s$ be positive integers with $x+y=n$ sufficiently large, $\lfloor n/2 \rfloor-s \le  x \le \lceil n/2\rceil +s$ and $s<n/10$.  Then
$$6{y \choose 4}+ 6(x-3){y \choose 3} \ge c(n,\bF) -(s+3)n^3.$$
\end{lemma}
\proof Define $f(y)=6{y \choose 4}+6(n-y-3){y \choose 3}$ and $a=\lfloor n/2 \rfloor$.  Then $c(n, \bF)=f(a)$.  We first observe that $f(y)$ is increasing for $1<y<n-4$. Indeed,
$$f(y+1)-f(y)=6\left({y \choose 3}+(n-y-4){y+1 \choose 3}-(n-y-3){y \choose 3}\right)=6{y \choose 2}(n-y-4)$$
and the condition on
$y$ shows that this is positive.
The condition  $\lfloor n/2 \rfloor -s \le x\le \lceil n/2\rceil +s$ implies that $1<a-s\le y \le a+1+s<n-4$ and so $f(y)\ge f(a-s)$. Therefore
\begin{align}
c(n, \bF)-f(y)
&\le f(a)-f(a-s)\notag \\
&=6\left({a \choose 4}+(n-a-3){a \choose 3}-{a-s \choose 4}-(n-a+s-3){a-s \choose 3}\right)\notag \\
&\le6\left(\frac{a^4}{4!}-\frac{(a-s-3)^4}{4!}+(n-a-3)\frac{a^3}{6}-(n-a+s-3)
\frac{(a-s-2)^3}{6}\right)\notag \\
&\le 6\left(\frac{a^4}{4!}-\frac{a^4-4a^3s-12a^3}{4!}+\frac{(n-a-3)}{6}(a^3-
(a^3-3a^2s-6a^2)\right)\label{middle} \\
&= 6\left(\frac{4a^3s+12a^3}{4!}+\frac{(n-a-3)(3a^2s+6a^2)}{6}\right)\notag\\
&<a^3s+3a^3+3a^2sn+6a^2n \notag \\
&<(s+3)n^3. \notag
\end{align}
Note that (\ref{middle}) follows from the inequalities
$(a-b)^4>a^4-4a^3b$ and $(a-b)^3>a^3-3a^2b$ which hold for $0<b<3a/2$; since $s<n/10$ we have $0<s+2<3a/2-1$.
This completes the proof of the Lemma.
\qed

We will need the following stability result proved independently by Keevash-Sudakov \cite{KSFano}
and F\"uredi-Simonovits \cite{FS}.

\begin{theorem} {\bf (${\bf F}$ Stability \cite{FS, KSFano})} \label{Fanostability}
Let $\cH$ be a 3-graph with $n$ vertices and $p^3(n)-o(n^3)$ edges
that contains no copy of ${\bf F}$. Then there is a partition of the
vertex set of $\cH$ into $X \cup Y$ so that the number of edges that
are within $X$ or within $Y$ is $o(n^3)$. In other words, $\cH$ can
be obtained from $P^3(n)$ by adding and deleting a set of $o(n^3)$
edges.
\end{theorem}

{\bf Remark.} The $o(1)$ notation above should be interpreted in the obvious way, namely $\forall \beta,  \exists \gamma, n_0$ such that if $n>n_0$ and $|\cH| >p^3(n)-\gamma n^3$, then $\cH=P^3(n) \pm \beta n^3$ edges.  We will
not explicitly mention the role of $\beta, \gamma$ when we use the result, but it should be obvious from the context. A similar comment applies for all applications of Theorem \ref{removal}.

 \bigskip

\noindent{\bf Proof of Theorem \ref{Fanoq}.}  Let $0< \delta \ll \e \ll 1$.
Write $o_{\delta}(1)$ for any function that approaches zero as
$\delta$ approaches zero and moreover, $o_{\delta}(1) \ll \e$. We emphasize that $\e$ is an absolute constant. Let
$n$ be sufficiently large and let $\cH$ be an $n$ vertex 3-graph
with $p^3(n)+q$ edges with $q <\e n^2$.  Write $\#{\bf F}$ for the number of copies of
${\bf F}$ in $\cH$.

If $\#{\bf F}\ge n^6$, then since $c(n, \bF)<n^4$, we have $\#\bF >\e n^2c(n,\bF)\ge qc(n, \bF)$ and we are done so assume that
$\#{\bf F}<n^6=(1/n) n^7$.  Since $n$ is sufficiently large, by the Removal lemma there is a set of at
most $\delta n^3$ edges of $\cH$ whose removal results in a
3-graph ${\cH}'$ with no copies of ${\bf F}$. Since
$|{\cH}'|>p^3(n)-\delta n^3$, by Theorem \ref{Fanostability},
we conclude that there is a bipartition of ${\cH}'$ (and also of
$\cH$) such that the number of edges contained entirely within a
part is $o_{\delta}(n^3)$. Now pick a bipartition $X \cup Y$ of
$\cH$ that maximizes $e(X,Y)$, the number of edges that intersect
both parts. We know that $e(X,Y)\ge p^3(n)-o_{\delta}(n^3)$, and an
easy calculation also shows that each of $X, Y$ has size
$n/2\pm o_{\delta}(n)$.

Let $B$ be the set of edges of $\cH$ that lie entirely within $X$ or entirely within $Y$ and let $G=\cH-B$. Let $M$ be the set of
triples which intersect both parts that are not edges of $\cH$. Then
$G \cup M$ is  bipartite so it has at most $p^3(n)$ triples. Consequently,
$$q+|M| \le |B| \le \odel(n^3).$$
Also, $|\cH|=|G|+|B|$ so we may suppose that $|G|=p^3(n)-t$ and $|B|=q+t$ for some $t \ge 0$.
For an edge $e \in B$, let $\bF(e)$ be the number of copies of $\bF$ in $\cH$ containing the unique edge $e$ from $B$.

If $t=0$, then $G \cong P^3(n)$ and $\bF(e)\ge c(n, \bF)$ for every $e \in B$ (by definition of $c(n, \bF)$) so we immediately obtain $\#\bF \ge qc(n, \bF)$. If $q>q(n, \bF)$ and $\bF(e)=c(n, \bF)$ for every $e \in B$, then there are two edges $e,e' \in B$ such that $|e \cap e'|=1$. To see this when $n$ is even, observe that if no two such edges exist, then every two edges of $B$ within $X$ intersect in zero or two points, and the same holds for the edges of $B$ within $Y$. The maximum number of edges that one can add to $P^3(n)$ with this property is $q(n, \bF)$, as every component is either a subset of $K_4^3$ or a sunflower with core of size two.  For $n$ odd we can only have edges in the larger part and again the same argument applies.
 
 We deduce that the number of copies of $\bF$ containing $e$ or $e'$ is at least $\bF(e)+\bF(e') +\bF(e,e')$ where
$\bF(e,e')$ is the  number of copies of $\bF$ in $\cH$ containing both $e$ and $e'$.
It is easy to see that $\bF(e, e') \ge 1$ (in fact, we have $\bF(e,e') \ge 2{\lfloor n/2 \rfloor\choose 2}$).

 We may therefore assume that $t\ge 1$ and we will now show that $\#\bF >qc(n, \bF)$. Partition $B=B_1 \cup B_2$, where
$$B_1=\{e \in B: \bF(e)>(1-\e)c(n,\bF)\}.$$
A potential copy of $\bF$ is a copy of $\bF$ in $G \cup M \cup B$ that uses exactly one edge of $B$.

\medskip

{\bf Claim 1.} $|B_1| \ge (1-\e)|B|$

Proof of Claim. Suppose to the contrary that $|B_2|\ge \e|B|$.
Pick $e=uvw \in B_2$.
Write $B_2=B_{XXX} \cup B_{YYY}$,
where the subscripts have the obvious meaning.  Assume by symmetry that $e \in B_{XXX}$.
For each  $Y'=\{y_1, \ldots, y_4\} \in {Y \choose 4}$, we can form a copy of ${\bf F}$ as follows:
Partition the six pairs of $Y'$ into three perfect matchings $L_u=\{e_u, e_u'\}, L_v=\{e_v, e_v'\}, L_w=\{e_w, e_w'\}$ and for each $x\in e$, add the two triples $x \cup e_x$ and $x \cup e_x'$.  There are six ways to choose the matchings $L_u, L_v, L_w$, so each choice of $Y'$ gives six potential copies of ${\bf F}$ containing $e$. Altogether we obtain $6{|Y| \choose 4}$ potential copies of ${\bf F}$. The only other way to form a copy of ${\bf F}$ using $e$ is to pick four points $a,b,c,d$ with $a \in X-e$ and $\{b,c,d\} \in {Y\choose 3}$. Then proceeding as before, we obtain $6(|X|-3){|Y| \choose 3}$ copies of ${\bf F}$.  This gives  a total of $(1-\od1)c(n,{\bf F})$ potential copies of $\bF$ containing $e$.
 At
least $(\e/2)c(n,{\bf F})$ of these potential copies of ${\bf F}$  have a
triple from $M$, for otherwise
$$\bF(e) \ge(1-o_{\delta}(1)-\e/2)c(n,{\bf F})>(1-\e)c(n,{\bf F})$$ which contradicts the definition of $B_2$.   The triple from $M$ referenced above lies in at most $2(|X||Y|+{|Y| \choose 2})<n^2$
copies of ${\bf F}$,  so the number of triples in $M$ counted here is at least
$${(\e/2)c(n,{\bf F})\over n^2}>(\e/30)n^2.$$
At least a third of these triples from $M$ are incident with the same vertex of $e$, so we conclude that
there exists $x \in e$ such that
$d_M(x)>(\e/100)n^2$.  Let $V=X \cup Y$ and let
$$A=\{v \in V: d_M(v)>(\e/100)n^2\}.$$
We have argued above that every $e \in B_2$ has a vertex in $A$.
Consequently,
$$3\sum_{v \in A}d_{B_2}(v) \ge 3|B_2|\ge 3\e|B|> 3\e|M|\ge \e\sum_{v
\in A} d_M(v) >\e|A|(\e/100)n^2,$$ and there exists a vertex $u \in  A$  such that
$d_{B_2}(u) \ge (\e^2/300)n^2$. Assume wlog that $u \in X$ so that
$d_{B_{XXX}}(u)\ge (\e^2/300)n^2$.

Let $\cH_{XYY}$ be the set of edges in $\cH$ with exactly one point in $X$. We may assume that $d_{H_{XYY}}(u)\ge d_{B_{XXX}}(u)$, for otherwise we may move $u$ to $Y$ and increase $e(X,Y)$, thereby contradicting the choice of $X, Y$. Consider $$e=uvw,\quad f=uy_1y_2,\quad f'=uy_1'y_2',$$ with $e \in B_{XXX}$ and $f, f' \in \cH_{XYY}$, $f \cap f'=\{u\}$. The number of choices
of $(e,\{f,f'\})$ is at least
$$d_{B_{XXX}}(u)\times \left({d_{\cH_{XYY}}(u) \choose 2} -n^3\right)>\e_1 n^6$$
where $\e_1=\e^6/10^{10}$.
 If for at least half of the choices  of $(e, \{f, f'\})$, these three edges span at least one copy of ${\bf F}$, then
 $\#{\bf F}>(\e_1/2)n^6>qc(n,\bF)$, a contradiction. So for at least half of the choices of $(e, \{f, f'\})$ above, $e \cup f \cup f'$ do not span a copy of ${\bf F}$. This implies that at least one of the triples  $xyy' \in M$ where $x \in e-\{u\}, y \in f-\{u\}, y' \in f' -\{u\}$.
  Since each such triple of $M$ is counted at most
$|X||Y|^2<n^3$ times, we obtain the contradiction
$(\e_1/2)n^6/n^3<|M|=o_{\delta}(n^3)$. This concludes the proof of the Claim. \qed
\bigskip

If $t\ge 4\e q$, then counting copies of $\bF$ from edges of $B_1$ and using Claim 1 we get
\begin{align}
\#\bF \ge \sum_{e \in B_1}(1-\e)c(n,\bF)&\ge |B_1|(1-\e)c(n,\bF) \notag \\
&\ge (1-\e)^2|B|c(n, \bF)\notag \\
&>(1-2\e)(q+t)c(n,\bF)\notag \\
&\ge (q+2\e q-8\e^2q)c(n,\bF)
\ge qc(n,\bF)\notag
\end{align}
and we are done. So we may assume that $t <4\e q<4\e^2n^2$. Let $x=|X|, y=|Y|$ and $s=\left\lceil \sqrt{2t/(n-2)}\right\rceil$.

{\bf Claim 2.} $\lfloor n/2 \rfloor-s \le  x \le \lceil n/2\rceil +s$ and if $t<(n-2)/2$, then $\lfloor n/2 \rfloor-s <  x < \lceil n/2\rceil +s$.

Proof of Claim.  We know that
$$p^3(n)-t=|G| \le {x \choose 2}y+{y \choose 2}x.$$
Now the Claim follows immediately from Lemma \ref{lemma1}.

Observe that $|M| \le t$ for otherwise $|G\cup M|>p^3(n)$ which is impossible.
Pick $e \in B$ and assume wlog that $e \subset X$. Since $t>0$, we have
$1 \le s\le \sqrt{2t/(n-2)}+1<n/10$.
The number of potential copies of $\bF$ containing $e$, denoted $pot\bF(e)$, is   $6{y \choose 4}+6(x-3){y \choose 3}$. Now Claim 2,
 Lemma \ref{lemma2} and $s\ge 1$ imply that
 $$pot\bF(e) \ge c(n,\bF)-(s+3)n^3\ge c(n, \bF)-4sn^3.$$ Not all of these copies of $\bF$ are in $\cH$, in fact, a triple from $M$ lies in at most $2n^2$ potential copies counted above (we pick either two more vertices in $Y$ or one in each of $Y$ and $X$, and there are two ways to complete a potential copy of \bF\, containing $e$). We conclude that
\begin{equation} \label{fe}
\bF(e) \ge pot\bF(e)-2n^2|M|\ge c(n, \bF)-4sn^3-2n^2|M| \ge c(n, \bF)-4sn^3-2tn^2.\end{equation}

Suppose first that $t< (n-2)/2$.  Then Claim 2 gives $\lfloor n/2 \rfloor-s <  x < \lceil n/2\rceil +s$.  Since $s=1$
and  $x$ is an integer, $|x-n/2|< 1$. The definition of $c(n,\bF)$ now yields
$$pot\bF(e) \ge \min \left\{
6{y \choose 4}+6(x-3){y \choose 3}: x \in \{\lfloor n/2 \rfloor, \lceil n/2 \rceil\}
\right\} \ge c(n, \bF).$$  Consequently, we can refine the bound in (\ref{fe}) to
$$\bF(e) \ge c(n, \bF)-2tn^2.$$  Altogether,
$$\#\bF \ge \sum_{e \in B}\bF(e) \ge (q+t)(c(n, \bF)-2tn^2)=qc(n,\bF)+tc(n,\bF)-2qtn^2-2t^2n^2.$$
Let us recall that $q \le \e n^2$ and $0<t<4\e q$. Then $2qtn^2<2\e tn^4<(t/2)c(n, \bF)$ and
$2t^2n^2=2t(tn^2)<(8\e q)tn^2<8\e^2tn^4<(t/2)c(n,\bF)$.
Consequently, $\#\bF > qc(n, \bF)$ as required.

Next we suppose that $t\ge (n-2)/2>n/4$.  This implies that $s\le  \sqrt{2t/(n-2)}+1\le 4\sqrt{t/n}$
and $\sqrt t \le 2t/\sqrt{n}$.  Therefore
$$4qsn^3<16qn^3\sqrt{t/n}=16q\sqrt t n^{2.5}\le 32qtn^2\le 32\e tn^4<(t/5)c(n, \bF).$$
So we again use (\ref{fe}) to deduce that $\#\bF$ is at least
$$ \sum_{e \in B}\bF(e) \ge (q+t)(c(n, \bF)-4sn^3-2tn^2)
\ge qc(n, \bF)+tc(n, \bF)-4qsn^3-2qtn^2-4tsn^3-2t^2n^2.$$
As $t<q<\e n^2$ we have the bounds
$$2qtn^2<(t/5)c(n,\bF), \quad 4tsn^3<4qsn^3<(t/5)c(n, \bF), \quad 2t^2n^2<(t/5)c(n, \bF).$$
This shows that $\#\bF>qc(n, \bF)$ and completes the proof of the theorem. \qed

\medskip

We end this section by proving that this result is asymptotically sharp.
\bigskip

\noindent{\bf Proof of Proposition \ref{pasch}.}
Let $0< \delta \ll \e$.
 Consider the following construction: Add a collection of $q$
edges to $P^3(n)$ within the part of size $\lceil n/2 \rceil$  such that the
following two conditions hold.

(1) every two added edges have at most
one point in common and

(2) the added edges do not form a Pasch
configuration, which is the six vertex 3-graph obtained from ${\bf F}$
by deleting a vertex.

It is well-known that such triple systems
exist of size $\delta n^2$ (in fact such Steiner triple systems also
exist \cite{GGW}).  Each new edge lies in at most
$c(n,\bF)$ copies of ${\bf F}$ that contain a unique new edge.  Now suppose that two of these
new edges, say $e,e'$ lie in a copy $C$ of ${\bf F}$.  Then there are at most
$n^2$ choices for the remaining two vertices of $C$. So the number of copies
of ${\bf F}$ containing two new edges   is at most $q^2n^2\le\delta q n^4 < \e qc(n,\bF)$. There are no copies of ${\bf F}$ using three new edges since three edges of ${\bf F}$ either span seven vertices or form a Pasch configuration. In either case we would have a Pasch configuration among the added edges.
Consequently, the number of copies of ${\bf F}$ is at most
$q(1+\e)c(n,\bF)$. \qed
\section{Counting  $F_5$'s}

Theorem \ref{F5q} follows from the following result.  Recall that $c(n, F_5)=(3+o(1))(n/3)^2$.

\begin{theorem} \label{F5}
For every $\e>0$ there exists $\delta>0$ and $n_0$ such that the
following holds for $n>n_0$. Every $n$-vertex 3-graph with
$t^3(n)+1$ edges contains either

$\bullet$ an edge that lies in at least $(3-\e)(n/3)^2$ copies of
$F_5$, or

$\bullet$ at least $\delta n^3$ copies of $F_5$.
\end{theorem}

\bigskip
\noindent{\bf Proof of Theorem \ref{F5q}.} Remove $q-1$ edges from
$\cH$ and apply Theorem \ref{F5}. If we find $\delta n^3$ copies of
$F_5$, then since $q<\delta n$, the number of copies is much
larger than $q(1-\e)c(n, F_5)$ and we are done. Consequently, we find
an edge $e_1$ in at least $(3-\e)(n/3)^2>(1-\e)c(n, F_5)$ copies of $F_5$. Now remove $q-2$ edges from $\cH-e_1$ and repeat this argument to obtain $e_2$.  In this way we obtain edges $e_1, \ldots, e_q$ as required.

 Sharpness
follows by adding a 3-partite triple system to one of the parts of
$T^3(n)$. It is easy to see that each added edge lies in
$c(n, F_5)-O(1)$ copies of $F_5$ and no copy of $F_5$ contains two of
the new edges. Consequently, the copies of $F_5$ are counted exactly
once. \qed

We will need the following stability theorem for $F_5$ proved by Keevash and the first author \cite{KM}.

\begin{theorem} {\bf($F_5$ Stability \cite{KM})}   \label{F5stability}
Let $\cH$ be a 3-graph with $n$ vertices and $t^3(n)-o(n^3)$ edges
that contains no copy of $F_5$. Then there is a partition of the
vertex set of $\cH$ into three parts so that the number of edges
with at least two vertices in some part is $o(n^3)$. In other words,
$\cH$ can be obtained from $T^3(n)$ by adding and deleting a set of
$o(n^3)$ edges.
\end{theorem}

 \bigskip

\noindent{\bf Proof of Theorem \ref{F5}.}
 Given $\e$ let $0<\delta \ll \e$.
Write $o_{\delta}(1)$ for a function that approaches zero as
$\delta$ approaches zero and moreover, $o_{\delta}(1) \ll \e$ for the set of functions used in this proof. Let
$n$ be sufficiently large and let $\cH$ be an $n$ vertex 3-graph
with $t^3(n)+1$ edges.  Write $\#F_5$ for the number of copies of
$F_5$ in $\cH$.

We first argue that we may assume that $\cH$ has minimum degree at
least  $d=(2/9)(1-\delta_1){n\choose 2}$, where
$\delta_1=\delta^{1/4}$. Indeed, if this is not the case, then
remove a vertex of degree less than $d$ to form the 3-graph $\cH_1$
with $n-1$ vertices.  Continue removing a vertex of degree less than
$d$  if such a vertex exists. If we could continue this process for
$\delta_2 n$ steps, where $\delta_2=\delta^{1/2}$, then the
resulting 3-graph $\cH'$ has $(1-\delta_2)n$ vertices and number of
edges at least
\begin{align}
{2\over 9}(1-\delta/2){n \choose 3}-(\delta_2n){2\over 9}(1-\delta_1){n \choose 2}
&\ge
{2\over 9}(1-\delta-3\delta_2(1-\delta_1)){n \choose 3} \notag \\
&>
{2\over 9}(1+\delta)(1-\delta_2)^3{n \choose 3} \notag \\
&>{2\over
9}(1+\delta){(1-\delta_2)n \choose 3}.\notag \end{align}
 By the result of
Keevash-Mubayi \cite{KM} and Erd\H os-Simonovits supersaturation
we conclude that $\cH$ has at least $\delta'n^5$ copies of $F^5$
(for some fixed $\delta'>0$) and we are done.  So we may assume that
this process of removing vertices of degree less than $d$ terminates
in fewer than $\delta_2 n$ steps, and when it terminates we are left
with a 3-graph $\cH'$ on $n'>(1-\delta_2)n$ vertices and  minimum
degree at least $d$.

Now suppose that we could prove that there is an edge of $\cH'$ that
lies in at least $(3-\e/2)(n'/3)^2$ copies of $F_5$. Since $\delta
\ll \e$, this is greater than $(3-\e)(n/3)^2$ and we are done. If on
the other hand $\cH'$ contains at least $2\delta n'^3$ copies of
$F_5$, then again this is at least $\delta n^3$ and we are done. So
if we could prove the result for $\cH'$ with $2\delta, \e/2$, then
we could prove the result for $\cH$ (with $\delta, \e$).
Consequently, we may assume that $\cH$ has minimum degree at least
$(2/9-o_{\delta}(1)){n \choose 2}=(1-\odel(1))(n/3)^2$.

If $\#F_5\ge \delta n^5$, then we are done so assume that
$\#F_5<\delta n^4$.  Then by the Removal lemma, there is a set of at
most $o_{\delta}(n^3)$ edges of $\cH$ whose removal results in a
3-graph ${\cH}'$ with no copies of $F_5$. Since
$|{\cH}'|>t^3(n)-o_{\delta}(n^3)$, by Theorem \ref{F5stability}, we
conclude that there is a 3-partition of ${\cH}'$ (and also of $\cH$)
such that the number of edges with at least two points in a part is
$o_{\delta}(n^3)$.   Now pick a partition $X \cup Y \cup Z$ of $\cH$
that maximizes $e(X,Y,Z)=\cH \cap (X\times Y \times Z)$. We know
that $e(X,Y,Z)\ge t^3(n)-o_{\delta}(n^3)$, and an easy calculation
also shows that each of $X, Y, Z$ has size $n/3+o_{\delta}(n)$.

Let $B={\cH} - (X \times Y \times Z)$ be the set of edges of $\cH$
that have at least two points in one of the partition classes and set $G=\cH -B$. Let
$M=(X\times Y \times Z) -\cH$ be the set of triples with one point
in each of $X,Y,Z$ that are not edges of $\cH$. Then $G \cup M=(\cH-B) \cup M$
is 3-partite so it has at most $t^3(n)$ triples.  Since $|\cH|=t^3(n)+1$, we conclude that
$$0\le |M|<|B| =o_{\delta}(n^3).$$

\noindent{\bf Claim.} For every vertex $v$ of $\cH$ we have $d_M(v)<\e'(n/3)^2$ for $\e'=\e/10^6$.

{\bf Proof of Claim.} Suppose for contradiction that $d_M(v)\ge \e'(n/3)^2$ for some vertex $v$.  Then
$$(1-\od1)(n/3)^2\le d_{\cH}(v)=d_G(v)+d_B(v)\le (1+\od1)(n/3)^2-\e'(n/3)^2+d_B(v).$$
We conclude that $d_B(v)\ge (\e'-\od1)(n/3)^2 > (\e'/2)(n/3)^2$.  Assume wlog that $v \in X$.

Case 1: $d_{B_{XXX}}(v)>(\e'/10)(n/3)^2$.
Suppose that $e=uvw$ satisfies $v \in e \in B_{XXX}$ and $(y, z) \in Y
\times Z$. The number of such choices for $(e, (y,z))$ is at least $d_{B_{XXX}}(v)|Y||Z|>(\e'/20)(n/3)^4$.
If for at least half of these choices  $e \cup \{y,z\}$
forms a copy of $F_5$ via the edges $e, uyz, wyz$  then we have $\#F_5>(\e'/40)(n/3)^4>\delta
n^3$, a contradiction. So for at least half of the choices of
$(e,(y,z))$ above, $xyz \not\in \cH$ for some $x \in \{u,w\}$  (i.e.
$xyz \in M$). Since each such triple of $M$ is counted at most
$|X|<n$ times (as $v$ is fixed), we obtain the contradiction
$(\e'/40n)(n/3)^4<|M|=o_{\delta}(n^3)$.  This concludes the proof in
this case.
\medskip

Case 2: $d_{B_{XXY}}(v)>(\e'/10)(n/3)^2$ or $d_{B_{XXZ}}(v)>(\e'/10)(n/3)^2$.
Assume by symmetry that $d_{B_{XXY}}(v)>(\e'/10)(n/3)^2$.
We may assume that $d_G(v) \ge d_{B_{XXY}}(v)$ for otherwise we can move $v$ to $Z$ and contradict the choice of the partition.
  Suppose that $e=uvw$ satisfies $v \in e \in B_{XXY}$ with $u \in X, w \in Y$. Let $(y, z) \in (Y-\{w\})
\times Z$ be such that $vyz \in \cH$. The number of such choices for $(e, (y,z))$ is at least $d_{B_{XXY}}(v)(d_G(v)-|Z|)>(\e'/11)^2(n/3)^4$.
If for at least half of these choices  $e \cup \{y,z\}$
forms a copy of $F_5$ via the triples $e, uyz, vyz$  then we have $\#F_5>(\e'/20)^2(n/3)^4>\delta
n^3$, a contradiction. So for at least half of the choices of
$(e,(y,z))$ above, $uyz \not\in \cH$ (i.e.
$uyz \in M$). Since each such triple of $M$ is counted at most
$|Y|<n$ times (as $v$ is fixed), we obtain the contradiction
$(\e'/20)^2(n/3)^4/n<|M|=o_{\delta}(n^3)$.  This concludes the proof in
this case.
\medskip

Case 3: $d_{B_{XYY}}(v)>(\e'/10)(n/3)^2$ or $d_{B_{XZZ}}(v)>(\e'/10)(n/3)^2$.
Assume by symmetry that $d_{B_{XYY}}(v)>(\e'/10)(n/3)^2$.
  Suppose that $e=uvw$ satisfies $v \in e \in B_{XXY}$ with $u,w \in Y$. Pick $(x, z) \in (X-\{v\})
\times Z$. The number of such choices for $(e, (x,z))$ is at least $d_{B_{XYY}}(v)(|X|-1)|Z|>(\e'/11)^2(n/3)^4$.
If for at least half of these choices  $e \cup \{x,z\}$
forms a copy of $F_5$ via the triples $xzu, xzw, e$  then we have $\#F_5>(\e'/20)^2(n/3)^4>\delta
n^3$, a contradiction. So for at least $(\e'/20)^2(n/3)^4$ of the choices of
$(e,(x,z))$ above, $xyz \not\in \cH$ for some $y \in \{u,w\}$ (i.e.
$xyz \in M$). For at least half of these choices, we may assume that $y=u$.
Since each such triple of $M$ is counted at most
$|Y|<n$ times (as $v$ is fixed), we obtain the contradiction
$(\e'/20)^2(n/3)^4/2n<|M|=o_{\delta}(n^3)$.  This concludes the proof of the Claim.

Let $B_1=B_{XXX} \cup B_{YYY} \cup B_{ZZZ} \subset B$, where the
subscripts have the obvious meaning ($B_{XXX}$ is the set of edges
in $B$ with three points in $X$ etc.), and let $B_2=B-B_1$, so $B_2$
consists of those edges of $\cH$ that have two points in one part
and one point in some other part.

Suppose that $e=uvw \in B_{XXX}$. For each  $(y,z) \in Y \times Z$
the points $u,v,w,y,z$  form a potential copy of $F_5$ via $e$ and two triples involving $y,z$.
For at least $(\e/2)(n/3)^2$ of these potential copies, $xyz \in M$ for $x \in e$, otherwise $e$ lies in
$(3-\od1-\e/2)(n/3)^2>(3-\e)(n/3)^2$ copies of $F_5$ and we are done.
 Each such triple of $M$ is counted at most twice, hence the number of triples intersecting $e$ is at least
$(\e/4)(n/3)^2$, and at least a third of these triples contain the same vertex $x \in e$. We conclude that
$d_M(x)>(\e/12)(n/3)^3\ge \e'(n/3)^2$ which contradicts the Claim.  The argument above works for any $e \in B_1$, so we have shown that $B_1=\emptyset$.

Let $e=uvw \in B_2=B$, where $u,v$ are in the same part, say $X$,
and $w$ is in another part, say $Y$. For each $(y,z) \in (Y-\{w\})
\times Z$, there are three types of potential copies of $F_5$ with
vertices $u,v,w,y,z$:

Type 1: $uyz, vyz, e$

Type 2: $uwz, e, vyz$ or $vwz, e, uyz$

The number of Type $i$ potential copies of $F_5$ is
$(|Y|-1)|Z|=(1-o_{\delta}(1))(n/3)^2$.  We may assume that the number of Type $1$  (real, not potential) copies of $F_5$
is at most $(1-\e/3)(n/3)^2$,  or that the number of Type $2$  (real, not potential) copies of $F_5$
is at most $(2-2\e/3)(n/3)^2$.
Otherwise $e$ lies in at least
$(3-\e)(n/3)^2$ copies of $F_5$ and we are done.

\medskip

Suppose that
the number of Type 1 copies of $F_5$ is at most
$(1-\e/3)(n/3)^2$. The number of pairs $(y,z) \in (Y-\{w\}) \times Z$ for
which either $uyz \in M$ or $vyz \in M$ is at least
$$(|Y|-1)|Z|-(1-\e/3)(n/3)^2 > (1-o_{\delta}(1)-1+\e/3)(n/3)^2>(\e/4)(n/3)^2.$$
Hence there exists $x \in \{u,v\}$ such that $xyz \in M$ for at
least $(\e/8)(n/3)^2$ pairs $(y,z) \in Y \times Z$. In other words,
$d_M(x)>(\e/8)(n/3)^2\ge \e'(n/3)^2$.  This contradicts the Claim.

We may therefore suppose that the number of Type 2 copies of $F_5$ is at most
$(2-2\e/3)(n/3)^2$. Assume by symmetry that there are at most $(1-\e/3)(n/3)^2$ Type 2 copies of the form
$uwz, e, vyz$. Arguing as above,  the number of pairs $(y,z) \in
Y \times Z$ for which either $uwz \in M$ or $vyz \in M$ is at least
$(\e/4)(n/3)^2$. If at least half of the time we have $vyz \in M$, then we obtain
$d_M(v)>(\e/8)(n/3)^2\ge \e'(n/3)^2$ and contradict the Claim. We therefore conclude that
for at least $(\e/8)(n/3)^2$ pairs $(y,z) \in Y \times Z$, we have $uwz \in M$. Consequently, the number of $z \in Z$ for which $uwz \in M$ is at least $(\e/10)(n/3)$.  We write this as $d_M(uw) \ge (\e/10)(n/3)$.

We have argued that for every edge $e=uvw \in B$ with $u,v$ in the same part and $w$ in a different part, either $d_M(uw) \ge \e n/30$ or $d_M(vw) \ge \e n/30$.  Form a bipartite graph with parts $B$ and $M$.  Let $e \in B$ be adjacent to $f \in M$ if $|e \cap f|=2$. We have shown above that each $e \in B$ has degree at least $\e n/30$.
Since $|B|>|M|$, we conclude that there exists $f \in M$ which is adjacent to at least $\e n/30$ different $e \in B$. Each of these $e \in B$ has two points in common point with $f$, so there is a pair of vertices $u,v$ in different parts of $\cH$ that lie is at least $\e n /90$ different $e \in B$.  Assume wlog that $u \in X, v \in Y$, and also that there are $x_i \in X$ for $1 \le i \le \e n /180$ such that $uvx_i \in B$ for each $i$.
For each $x_i$, consider $(y,z) \in (Y -\{v\})\times Z$ and triples $x_ivz, x_ivu, uyz$. The number of such choices for $(i,y,z)$ is at least $(\e n/200)(n/3)^2$.  If for at least half of these choices these three triples are edges of $\cH$, then we obtain $\#F_5 \ge (\e n/400)(n/3)^2>\delta n^3$ and we are done. So for at least half of these choices of $(i,y,z)$  we have either $x_ivz \in M$ or $uyz \in M$. Each such triple of $M$ is counted at most $n$
times so we obtain at least $(\e/400)(n/3)^2$ triples from $M$ incident to some vertex of $e$. At least one third of these triples are incident to the same vertex of $e$, so we obtain $x \in e$ with
$d_M(x)\ge (\e/1200)(n/3)^2\ge \e'(n/3)^2$.  The contradicts the Claim and completes the proof. \qed

\bigskip

\section{Counting $B_5$'s}

Theorem \ref{B5q} follows from the following result.  Recall that $c(n, B_5)=(2+o(1))(n/3)^2$.

\begin{theorem} \label{B5}
For every $\e>0$ there exists $\delta>0$ and $n_0$ such that the
following holds for $n>n_0$. Every $n$-vertex 3-graph with
$b^3(n)+1$ edges contains either

$\bullet$ an edge that lies in at least $(2-\e)(n/3)^2$ copies of
$B_5$, or

$\bullet$ at least $\delta n^4$ copies of $B_5$.
\end{theorem}

\bigskip
\noindent {\bf Proof of Theorem \ref{B5q}.} Remove $q-1$ edges from
$\cH$ and apply Theorem \ref{B5}. If we find $\delta n^4$ copies of
$B_5$, then since $q<\delta n^2$, the number of copies is much
larger than $(1-\e)c(n, B_5)$ and we are done.
Consequently, we find
an edge $e_1$ in at least $(2-\e)(n/3)^2>(1-\e)c(n, B_5)$ copies of $B_5$. Now remove $q-2$ edges from $\cH-e_1$ and repeat this argument to obtain $e_2$.  In this way we obtain edges $e_1, \ldots, e_q$ as required.

 Sharpness
follows by adding a partial Steiner triple system to $B^3(n)$ where
each added edge is entirely within $X$. In other words, we are
adding a collection of triples within $X$ such that every two have
at most one point in common.  It is easy to see that each added edge
lies in $c(n, B_5)-O(1)$ copies of $B_5$ and moreover, since these
edges have at most one common point, these copies are counted
exactly once. \qed

We will need the following stability theorem for $B_5$ proved by F\"uredi-Pikhurko-Simonovits \cite{FPS}.

\begin{theorem} {\bf ($B_5$ stability \cite{FPS})} \label{B5stability}
Let $\cH$ be a 3-graph with $n$ vertices and $b^3(n)-o(n^3)$ edges
that contains no copy of $B_5$. Then there is a partition of the
vertex set of $\cH$ into $X \cup Y$ so that the number of edges that
are not of the form $XXY$ is $o(n^3)$. In other words, $\cH$ can be
obtained from $B^3(n)$ by adding and deleting a set of $o(n^3)$
edges.
\end{theorem}

 \bigskip

 \noindent{ \bf Proof of Theorem \ref{B5}.}
 Given $\e$ let $0<\delta \ll \e$.
Write $o_{\delta}(1)$ for any function that approaches zero as
$\delta$ approaches zero and moreover, $o_{\delta}(1) \ll \e$. Let
$n$ be sufficiently large and let $\cH$ be an $n$-vertex 3-graph
with $b^3(n)+1$ edges.  Write $\#B_5$ for the number of copies of
$B_5$ in $\cH$.

We first argue that we may assume that $\cH$ has minimum degree at
least  $d=(4/9)(1-\delta_1){n\choose 2}$, where
$\delta_1=\delta^{1/4}$. Indeed, if this is not the case, then
remove a vertex of degree less than $d$ to form the 3-graph $\cH_1$
with $n-1$ vertices.  Continue removing a vertex of degree less than
$d$  if such a vertex exists. If we could continue this process for
$\delta_2 n$ steps, where $\delta_2=\delta^{1/2}$, then the
resulting 3-graph $\cH'$ has $(1-\delta_2)n$ vertices and number of
edges at least
$${4\over 9}(1-\delta-3\delta_2(1-\delta_1){n \choose 3}>{4\over
9}(1+\delta){(1-\delta_2)n \choose 3}.$$ By the result of
F\"uredi-Pikhurko-Simonovits \cite{FPS} and Erd\H os-Simonovits supersaturation
we conclude that $\cH$ has at least $\delta'n^5$ copies of $B^5$
(for some fixed $\delta'>0$) and we are done.  So we may assume that
this process of removing vertices of degree less than $d$ terminates
in at most $\delta_2 n$ steps, and when it terminates we are left
with a 3-graph $\cH'$ on $n'>(1-\delta_2)n$ vertices and  minimum
degree at least $d$.

Now suppose that we could prove that there is an edge of $\cH'$ that
lies in at least $(2-\e/2)(n'/3)^2$ copies of $B_5$. Since $\delta
\ll \e$, this is greater than $(2-\e)(n/3)^2$ and we are done. If on
the other hand $\cH'$ contains at least $2\delta n'^4$ copies of
$B_5$, then again this is at least $\delta n^4$ and we are done. So
if we could prove the result for $\cH'$ with $2\delta, \e/2$, then
we could prove the result for $\cH$ (with $\delta, \e$).
Consequently, we may assume that $\cH$ has minimum degree at least
$(4/9-o_{\delta}(1)){n \choose 2}$.

If $\#B_5\ge \delta n^4$, then we are done so assume that
$\#B_5<\delta n^4$.  Then by the Removal lemma, there is a set of at
most $o_{\delta}(n^3)$ edges of $\cH$ whose removal results in a
3-graph ${\cH}'$ with no copies of $B_5$. Since
$|{\cH}'|>b^3(n)-o_{\delta}(n^3)$, by Theorem \ref{B5stability}, we
conclude that there is a partition $X \cup Y$ of the vertex set of
${\cH}'$ (and also of $\cH$) such that the number of edges with $0,
1$, or $3$ points in $X$ is $o_{\delta}(n^3)$.   Now pick a
partition $X \cup Y$ of $\cH$ that maximizes $e(X,X,Y)$ the number
of edges with exactly two points in $X$.  We know that $e(X,X,Y)\ge
b^3(n)-o_{\delta}(n^3)$, and an easy calculation also shows that
$|X|=2n/3+o_{\delta}(n)$ and $|Y|=n/3+o_{\delta}(n)$.

Let  $B$ be the set of edges of $\cH$ that do not have exactly two
points in $X$. Let $M$ be the set of triples with exactly two points
in $X$ that are not edges of $\cH$ and let $G=\cH-B$ be the set of edges of $\cH$ with exactly two points in $X$. Then $\cH-B \cup M$ has a
$(2,1)$-partition $X \cup Y$, so it has at most $b^3(n)$ edges. We
conclude that
$$|M|<|B| =o_{\delta}(n^3).$$
In particular, $B \ne \emptyset$.  Partition $B=B_{XXX} \cup B_{XYY}
\cup B_{YYY}$, where $B_{X^iY^{3-i}}$ is the set of edges in $B$
with $i$ points in $X$ and $3-i$ points in $Y$.
\medskip

\noindent{\bf Claim 1.}  For every vertex $v$ of $\cH$ we have $d_{B_{XXX}}(v)<\e_1n^2$, where $\e_1=\e^2/10^6$.

\noindent
{\bf Proof of Claim 1.} Suppose for contradiction that
$d_{B_{XXX}}(v)>
\e_1n^2$ for some vertex $v$.  Let $B(v)$ be the set of edges in $B_{XXX}$ that
contain $v$, so $|B(v)|=d_{B_{XXX}}(v)$.
First  observe that $d_G(v) \ge d_{B_{XXX}}(v)$,  for
otherwise we can move $v$ to $Y$ and contradict the choice of the
partition $X, Y$.
Now for each $e=vab \in B(v)$ and $f=vxy \in \cH$ with $\{a,b,x\} \in {X-\{v\}\choose 3},
y \in Y$, consider the two triples $axy, bxy$. We see that $e, f,
axy, bxy$ forms a (potential) copy of $B_5$. For each $e$, the
number of $f$ is at least $d_{G}(v)-n \ge |B(v)|-n>|B(v)|/2$, since $f$ must omit
$a,b$ and there are at most $|Y|$ pairs containing either of them.
Hence the number of choices for $(e,f)$ is at least $|B(v)|^2/2$.
If for at least half of these choices  of $(e,f)$, we obtain a copy of
$B_5$ in $\cH$,
 then  $\#B_5>|B(v)|^2/4>\delta
n^4$, a contradiction. So for at least half of the choices of $(e,f)$
above, one of the triples $axy, bxy$ is in $M$. A given triple in
$M$ is counted  at most $|X|<n$ times, so we obtain the
contradiction $|B(v)|^2/(4n)<|M|=o_{\delta}(n^3)$.  This finishes the proof of the Claim.
\medskip

\noindent {\bf Case 1.} $|B_{XXX}|\ge |B|/3$.

For each $e=uvw \in B_{XXX}$, and $(x,y) \in (X-e) \times Y$, there
is a potential copy of $B_5$ consisting of vertices $u, v, w, x, y$
and edges $uxy, vxy, wxy, e$. This gives a total of
$(|X|-3)|Y|>(2-o_{\delta}(1))(n/3)^2$ potential copies of $B_5$. At
least $(2\e/3)(n/3)^2$ of these potential copies of $B_5$ have a
triple from $M$, for otherwise $e$ would lie in at least
$(2-o_{\delta}(1)-2\e/3)(n/3)^2>(2-\e)(n/3)^2$ copies of $B_5$ and
we are done. The triple from $M$ referenced above cannot be $e$
(since $e\in \cH$), and therefore lies in exactly one copy of $B_5$ that was counted above.
At least a third of these triples from $M$ are incident with the
same vertex of $e$, hence there exists $z \in e$ such that
$d_M(z)>(2\e/9)(n/3)^2$.

Let $V=X \cup Y$ and let
$$A=\{v \in V: d_M(v)>(2\e/9)(n/3)^2\}.$$
We have argued above that every $e \in B_{XXX}$ has a vertex in $A$.
Consequently,
$$9\sum_{v \in A}d_{B_{XXX}}(v)\ge
9|B_{XXX}|\ge 3|B|>3|M|\ge \sum_{v \in A} d_M(v)
>|A|(2\e/9)(n/3)^2,$$
and there exists a vertex $v \in X \cap A$  such that $d_{B_{XXX}}(v)>
(\e/50)(n/3)^2>\e_1n^2$. This contradicts Claim 1 and concludes
the proof in this case.

\medskip

\noindent {\bf Case 2.} $|B_{YYY}|\ge |B|/3$.

For each $e=uvw \in B_{YYY}$ and $x,x' \in X$, there is a potential
copy of $B_5$ consisting of vertices $u,v,w, x,x'$ and edges $xx'u,
xx'v, xx'w, e$. This gives a total of ${|X| \choose
2}>(2-o_{\delta}(1))(n/3)^2$ potential copies of $B_5$. At least
$(2\e/3)(n/3)^2$ of these potential copies of $B_5$ have a triple
from $M$, for otherwise $e$ would lie in at least
$(2-o_{\delta}(1)-2\e/3)(n/3)^2>(2-\e)(n/3)^2$ copies of $B_5$ and
we are done. The triple from $M$ referenced above cannot be $e$
(since $e\in \cH$), and therefore lies in exactly one copy of $B_5$.
At least a third of these triples from $M$ are incident with the
same vertex of $e$, hence there exists $z \in e$ such that
$d_M(z)>(2\e/9)(n/3)^2$.  As in Case 1,
let $V=X \cup Y$ and $A=\{v \in V: d_M(v)>(2\e/9)(n/3)^2\}$.
We have argued above that every $e \in B_{YYY}$ has a vertex in $A$.
Consequently,
$$9\sum_{v \in A}d_{B_{YYY}}(v)\ge
9|B_{YYY}|\ge 3|B|> 3|M|\ge \sum_{v \in A} d_M(v)
>|A|(2\e/9)(n/3)^2,$$
and there exists a vertex $v\in Y \cap A$  such that $d_{B_{YYY}}(v)>
(\e/50)(n/3)^2$.  Let $B(v)$ be the set of edges in $B_{YYY}$ that
contain $v$, so $|B(v)|=d_{B_{YYY}}(v)$.

Next we observe that $d_{B_{XYY}}(v) \le d_{G}(v)$ otherwise
we can move $v$ to $X$ and contradict the choice of the partition
$X, Y$.  We also recall that $\cH$ has minimum degree at least
$(4/9-o_{\delta}(1)){n \choose 2}$, so
$$d_{B_{YYY}}(v)+d_{G}(v)+d_{B_{XYY}}(v)\ge (4/9-o_{\delta}(1)){n \choose
2}.$$ Since $d_{B_{YYY}}(v)\le {|Y| \choose 2}<(1/9+o_{\delta}(1)){n
\choose 2}$, we conclude that
$$d_{G}(v)>{1\over 2}\left({4\over 9}-{1\over 9}-o_{\delta}(1)\right){n \choose 2}=
\left({1\over 6}-o_{\delta}(1)\right){n \choose 2}.$$

Now for each $e=vyy' \in B(v)$ and $f=vxx' \in G$  ($x,x' \in
X$), consider the two triples $xx'y, xx'y'$. We see that $e,f, xx'y,
xx'y'$ forms a potential copy of $B_5$. The number of choices of
$(e,f)$ above is at least $|B(v)| d_{G}(v)=d_{B_{YYY}}(v) d_G(v)$. If for at least
half of these choices of $(e,f)$, we obtain a copy of $B_5$ in $\cH$,
 then  $\#B_5>d_{B_{YYY}}(v) d_G(v)/2>\delta
n^4$, a contradiction. So for at least half of the choices of $(e,f)$
above, one of the triples $xx'y, xx'y'$ is in $M$. A given triple in
$M$ is counted  at most $|Y|<n$ times, so we obtain the
contradiction $d_{B_{YYY}}(v)d_{G}(v)/(2n)<|M|=o_{\delta}(n^3)$. This
concludes the proof in this case.

\medskip

\noindent {\bf Case 3.} $|B_{XYY}|\ge |B|/3$.

Let
$$B_1=\{e \in B_{XYY}: \hbox{there exists } v \in e \cap Y \hbox{
with } d_M(v)>\e (n/3)^2\}.$$

\noindent {\bf Subcase 3.1.} $|B_1| \ge |B_{XYY}|/2$. Let
$$A=\{v \in Y: d_M(v)>(\e/2)(n/3)^2\}.$$
By definition, every $e \in B_1$ has a vertex in $A$.  Therefore
$$18\sum_{v \in A}d_{B_{1}}(v)\ge 18|B_1|\ge 9|B_{XYY}|
\ge 3|B|> 3|M|\ge \sum_{v \in A} d_M(v)
>|A|(\e/2)(n/3)^2,$$
and there exists a vertex $v \in Y$  such that
$$d_{B_{XYY}}(v) \ge d_{B_{1}}(v)>
(\e/36)(n/3)^2.$$
Recall that $G$ is the set of edges of $\cH$ with
exactly two points in $X$.  Next observe that $d_{G}(v)\ge
d_{B_{XYY}}(v)$ for otherwise we can move $v$ to $X$ which increases
$e(X,X,Y)$ and  contradicts the choice of $X,Y$. It follows that
$d_{G}(v)>(\e/36)(n/3)^2$.

Now for each $e=uvw \in B_{XYY}$ and $f=xx'v \in G$ with
$\{u,x,x'\} \in {X \choose 3}$, and $w\in Y$, consider the two triples
$uwx, uwx'$. We see that $e,uwx, uwx', xx'v$ forms a potential copy
of $B_5$. The number of choices of $(e,f)$ above is at least
$d_{B_{XYY}}(v)\times(d_{G}(v)-|X|)>d_{B_{XYY}}(v)d_{G}(v)/2$.
If for at least half of these choices of $(e,f)$, we obtain a copy of
$B_5$ in $\cH$,
 then
 $$\#B_5>{d_{B_{XYY}}(v)d_{G}(v)\over 4}>{\e^2\over 10^5}\left({n\over 3}\right)^4>\delta
n^4,$$ a contradiction. So for at least half of the choices of $(e,f)$
above, one of the triples $uwx, uwx'$ is in $M$. A given triple in
$M$ is counted  at most $|X|<n$ times, so we obtain the
contradiction
$${\e^2\over 10^5}\left({n^3\over 3^4}\right)<{d_{B_{XYY}}(v)d_{G}(v)\over 4n}<|M|=o_{\delta}(n^3).$$ This
concludes the proof in this subcase.

\noindent {\bf Subcase 3.2.} $|B_1| < |B_{XYY}|/2$. So in this
subcase we have $|B_2| \ge |B_{XYY}|/2$, where
$$B_2=\{e \in B_{XYY}: \hbox{for every  } v \in e \cap Y \hbox{
we have } d_M(v)\le (\e/2) (n/3)^2\}.$$ Fix $e=uvw \in B_2$ with $u
\in X$ and $v,w \in Y$.

\noindent{\bf Claim 2.} There exist  sets $X_v, X_w \subset X$ such
that

$\bullet$ $xuv \in M$ for every $x \in X_v$ and $xuw \in M$ for
every $x \in X_w$ and

$\bullet$ $|X_v|>(\e/20)n$ and   $|X_w|>(\e/20)n$

\noindent{\bf Proof of Claim 2.} Let $X_v=\{x \in X: xuv \in M\}$. We
will show that $|X_v|\ge (\e/20)n$. The same argument will apply to
$X_w$.

 Suppose for contradiction that $|X_v|<(\e/20)n$.
Pick $x,x' \in X-X_v$ and consider $u,v,w,x,x'$. The triples $uvx,
uvx', e, xx'w$ form a potential copy of $B_5$.   Since $x,x' \in
X-X_v$, we have $uvx \in \cH$ and $uvx' \in \cH$. So if these four
edges do not form a copy of $B_5$ in $\cH$ then $xx'w \in M$. Since
$e \in B_2$, the number of pairs $\{x,x'\} \in {X \choose 2}$ such that $xx'w \in M$
is at most $(\e/2)(n/3)^2$. Consequently, the number of pairs $x,x'
\in X-X_v$ with $xx'w \in \cH$ is at least
\begin{align}
{|X-X_v| -1\choose 2}-
{\e\over 2}\left({n\over 3}\right)^2
&>{(1-o_{\delta}(1)-{3\e\over 40}){2n\over 3} \choose 2}-{\e\over 2}\left({n\over
3}\right)^2 \notag \\
&>\left(2\left(1-{\e\over 10}\right)^2 -{\e\over
2}\right)
 \left({n \over 3}\right)^2 \notag \\
 &=\left(2-{9\e \over 10}+{\e^2 \over
 50}\right)\left({n \over 3}\right)^2 \notag \\
&> (2-\e)\left({n \over 3}\right)^2\notag.\end{align}

This gives us the required number of copies of $B_5$ containing the
edge $e$ and concludes the proof of the Claim.

\medskip

For each edge $e=uvw \in B_2$ with $u \in X, v,w \in Y$, Claim 2
shows that are at least $(\e/20)n$ triples of the form $xuv \in M$.
Form the bipartite graph with parts $B_2$ and $M$, where $uvw \in B_2$ is
adjacent to all such $xuv \in M$. Then since every vertex of $B_2$ has
degree at least $(\e/20)n$, and $|B_2| \ge |B_{XYY}|/2 \ge |B|/6 >
|M|/6$, we conclude that there exists $xuv \in M$ (with $v \in Y$) which is adjacent
to at least $(\e/120)n$ edges in $B_2$. Each of these edges of $B_2$ contains $v$, and either $x$ or $u$, so we may assume by symmetry that at least half of them contain $u$. So we have
$uvw_i \in B_2$, where $u \in X$ and
$v, w_i \in Y$ for $i=1, \ldots, (\e/240)n$. For each $w_i$,
consider the set $X_{w_i}$ defined in Claim 2. We know that $x'
uw_i \in M$ for each $w_i$ and  $x' \in X_{w_i}$.  Since these
triples are distinct for distinct $w_i$ or distinct $x'$, we
conclude that $d_M(u)\ge (\e/240)n (\e/20)n=(\e^2/4800)n^2$.
Recalling the minimum degree condition on $\cH$, we have
$$(4/9 -o_{\delta}(1)){n \choose 2} \le d_{\cal
H}(u)=d_{G}(u)+d_{B}(u)\le (4/9 -o_{\delta}(1)){n \choose 2} -d_M(u)+d_B(u).
$$  We conclude that $d_B(u)\ge (\e^2/5000)n^2$.  By Claim 1 we know that
$d_{B_{XXX}}(u)<\e_1 n^2$ where $\e_1=\e^2/10^5$.  As $d_B(u)=d_{B_{XXX}}(u)+d_{B_{XYY}}(u)$, we obtain
$$d:=d_{B_{XYY}}(u)=d_B(u)-d_{B_{XXX}}(u)>(\e^2/5000)n^2-\e_1n^2\ge2\e_1n^2.$$ Say that $uyy' \in B_{XYY}$ is bad if
$$|\{ x \in X: xuy \in M \hbox{ or } xuy' \in M\}| >(1-\e_1)(2n/3).$$
Let
$$S=\left\{\{y, y'\} \in {|Y| \choose 2}: uyy' \hbox{ is bad }\right\}.$$
 Now suppose that $|S| \ge (0.9)d$. For each
 $e=uyy' \in B_{XYY}$ with $\{y,y'\} \in S$ there is a set $X_e
 \subset X$ with $|X_e| \ge (1-\e_1)(2n/3)$ such that $xuy \in M$ or $xuy' \in
 M$ for all $x \in X_e$. Each of these triples in $M$ is counted at
 most $|Y|$ times so we obtain
 \begin{equation} \label{m} d_M(u) \ge {|S|(1-\e_1)(2n/3) \over (1+\odel(1))n/3}
 =2(1-2\e_1)|S|\ge (1.8)(1-2\e_1)d>(1.7)d.\end{equation}
Again recalling the minimum degree condition on $\cH$, we have
\begin{align}
(4/9 -o_{\delta}(1)){n \choose 2} \le d_{\cal
H}(u)&=d_{B_{XXX}}(u)+d_{G}(u)+d_{B_{XYY}}(u) \notag \\
&=d_{B_{XXX}}(u)+\left((4/9+o_{\delta}(1)){n \choose
2}-d_{M}(u)\right)+d.\notag
\end{align} Using (\ref{m}) and $d>\e_1n^2$
we obtain $d_{B_{XXX}}(u)> (0.7)d-\odel(n^2)>\e_1n^2$.
 This contradicts Claim 1 and
concludes the proof if $|S| \ge (0.9)d$.

Next suppose that $|S| < (0.9)d$. So for at least $(0.1)d$ edges
$e=uyy' \in B_{XYY}$ we have a set $X_e \subset X$ such that
$$|X_e|\ge (\e_1-o_{\delta}(1))(2n/3) >(\e_1/3)n$$
and $uyx \in \cH$ for all $x \in X_e$ (also $uy'x \in \cH$ but we
wont use this).

Let $x,x'\in X_e$ and
consider the  triple $xx'y'$. We see that $e,uyx, uyx', xx'y'$ forms
a potential copy of $B_5$. The number of choices for $(e,\{x,x'\})$
above is at least
$$(0.1)d \times {(\e_1/3)n \choose 2}>(\e_1/5)n^2 \times (\e_1^2/20)n^2=(\e_1^3/100)n^4.$$ If for at
least half of these choices of $(e,\{x,x'\})$, we have $xx'y' \in
\cH$,
 then
 $$\#B_5>(\e_1^3/200)n^4>\delta
n^4,$$ a contradiction. So for at least half of the choices of
$(e,\{x,x'\})$ above, $xx'y' \in M$. A given triple $xx'y' \in M$ is
counted at most $|Y|<n/2$ times, so we obtain the contradiction
$${\e_1^3 \over 200}n^3 \le {(0.1)d \times {|X_e| \choose 2} \over n}<|M|=o_{\delta}(n^3).$$
This completes the proof of the subcase and the Theorem.
 \qed

\section{Counting Expansions of Cliques}

In this section we will prove Theorems \ref{Lq} and \ref{exactLq}.

\subsection{Asymptotic Counting}

Theorem \ref{Lq} follows from the following result.  Recall that $l_{r+1}={r+2 \choose 2}$ is the number of vertices of $L_{r+1}$ and $c(n, L_{r+1})=\Theta(n^{l_{r+1}-3})$.

\begin{theorem} \label{L1}
For every $\e>0$ there exists $\delta>0$ and $n_0$ such that the
following holds for $n>n_0$. Every $n$ vertex 3-graph with
$t^3_r(n)+1$ edges contains either

$\bullet$ at least $\delta n^{l_{r+1}-1}$ copies of $L_{r+1}$, or

$\bullet$ an edge that lies in at least $c(n, L_{r+1})$ copies of
$L_{r+1}$, or

$\bullet$ two edges that each  lie in at least $(1-\e)c(n, L_{r+1})$ copies of $L_{r+1}$ with none of these copies containing both edges.
\end{theorem}

\bigskip
\noindent {\bf Proof of Theorem \ref{Lq}.} Remove $q-1$ edges
from $\cH$ and apply Theorem \ref{L1}. If we find $\delta n^{l_{r+1}-1}$
copies of $L_{r+1}$, then since $q<\delta n^2$, the number of copies is
much larger than $q(1-\e)c(n, L_{r+1})$ and we are done.
Consequently, we find
an edge $e_1$ in at least $(1-\e)c(n, L_{r+1})$ copies of $L_{r+1}$. Now remove $q-2$ edges from $\cH-e_1$ and repeat this argument to obtain $e_2$.  In this way we obtain edges $e_1, \ldots, e_q$ as required.

Sharpness
follows by the following construction:  Take $T^3_r(n)$ with parts $V_1, \ldots, V_r$, pick any point $y \in V_2$, and add $q$ edges of the form $xx'y$ with $x,x' \in V_1$.  Each added edge lies in at most $(1+\e)c(n, L_{r+1})$ copies of $L_{r+1}$, and no two added edges lie in a common copy of $L_{r+1}$, since $L_{r+1}$ has the property that for every two edges $e,e'$ containing a common vertex $v$, there is another edge $f$ containing a point from each of $e-\{v\}$ and $e' -\{v\}$ and $v \not\in f$. Taking two edges containing $y$, we see that there is no edge that can play the role of $f$ above.
\qed
\medskip

We will need the following stability result proved by Pikhurko \cite{P} (see also \cite{M}).
\bigskip

\begin{theorem} {\bf ($L_{r+1}$ Stability \cite{P})} \label{Lstability}
Let $\cH$ be a 3-graph with $n$ vertices and $t^3_r(n)-o(n^3)$ edges
that contains no copy of $L_{r+1}$. Then there is a partition of the
vertex set of $\cH$ into $r$ parts so that the number of edges that
 intersect some part in at least two points is $o(n^3)$. In other words, $\cH$ can
be obtained from $T^3_r(n)$ by adding and deleting a set of $o(n^3)$
edges.
\end{theorem}

\bigskip

\noindent{\bf Proof of Theorem \ref{L1}.}
 Given $\e$ let $0<\delta \ll \e$.
Write $o_{\delta}(1)$ for any function that approaches zero as
$\delta$ approaches zero and moreover, $o_{\delta}(1) \ll \e$. Let
$n$ be sufficiently large and let $\cH$ be an $n$ vertex 3-graph
with $t^3_r(n)+1$ edges.  Write $\#L_{r+1}$ for the number of copies of
$L_{r+1}$ in $\cH$.

If $\#L_{r+1}\ge \delta n^{l_{r+1}-1}$, then we are done so assume that
$\#L_{r+1}<\delta n^{l_{r+1}-1}$.  Then by the Removal lemma, there is a set of at
most $o_{\delta}(n^3)$ edges of $\cH$ whose removal results in a
3-graph ${\cH}'$ with no copies of $L_{r+1}$. Since
$|{\cH}'|>t^3_r(n)-o_{\delta}(n^3)$, by Theorem \ref{Lstability},
we conclude that there is an $r$-partition $V_1 \cup \cdots \cup V_r$ of ${\cH}'$ (and also of
$\cH$) such that the number of edges that intersect some part in at least two points
 is $\odel(n^3)$.
  Now pick a partition $V_1 \cup \cdots \cup V_r$ of
$\cH$ that maximizes $h_1+2h_2+3h_3$, where $h_i$ is the number of edges of $\cH$ that intersect precisely $i$ of the parts. The partition guaranteed by Theorem \ref{Lstability} satisfies $h_1+2h_2=\odel (n^3)$, and
hence for this particular partition
$h_1+2h_2+3h_3\ge 3|\cH| -2(h_1+h_2)>3t^3_r(n)-\odel(n^3)$.
Since $h_1+2h_2+3h_3\le 3|\cH|-(h_1+h_2)$ we conclude that
 for the partition that maximizes $h_1+2h_2+3h_3$
we have $h_1+2h_2=\odel (n^3)$ and $h_3\ge t^3_r(n)-o_{\delta}(n^3)$. A
standard calculation also shows that for this partition each $V_i$ has size
$n/r\pm o_{\delta}(n)$.

Let $B=\cH - \prod_{i=1}^r V_i$, let  $G=\cH-B$  and $M=\prod_{i=1}^r V_i-G$. Then
$\cH-B \cup M$ is  $r$-partite so it has at most $t^3_r(n)$ edges. We
conclude that
$$|M|<|B| =o_{\delta}(n^3),$$
in particular $|B|\ge 1$.  We will now argue that we can improve this to $|B|\ge 2$.
We may suppose  that $n_i:=|V_i|$ satisfy
$n_1 \ge n_2 \ge \ldots \ge n_r$.
Pick $e_1 \in B$.  If $\cH-e_1 \cong T^3_r(n)$, then clearly $e_1$ lies in at least $c(n, L_{r+1})$ copies of $L_{r+1}$ and we are done.  So assume that $\cH-e_1 \not\cong T^3_r(n)$.
Suppose that $B \cap (\cH -e_1) =\emptyset$. Then either $n_r \ge n_1-1$ and
$$t^3_r(n)=|\cH-e_1| \le \left(\sum_{S \in {[r] \choose 3}}\prod_{i \in S}n_i\right)-1<t^3_r(n),$$
or $n_r<n_1-1$ and
$$t^3_r(n)=|\cH-e_1| \le \sum_{S \in {[r] \choose 3}}\prod_{i \in S}n_i<t^3_r(n).$$
In either case we have a contradiction, so we may assume that $B \cap (\cH -e_1) \ne\emptyset$.  In other words, there exists $e_2\ne e_1$ such that $e_2 \in B$ and therefore $|B| \ge 2$. We will now show that every $e \in B$ lies in at least $(1-\e)c(n, L_{r+1})$ copies of $L_{r+1}$ in $\cH$ and each copy uses a unique edge from $B$.

Let $e=xyz\in B$.  We may assume by symmetry that  $x,y \in V_1$. Pick $(v_2, \ldots, v_r) \in V_2 \times \cdots \times V_r$
with $v_i \ne z$ for all $i$. For every pair of distinct vertices $\{a,b\}$ with $a \in \{v_2, \ldots, v_r\}$ and $b \in\{x,y,v_2,\ldots, v_r\}$ (there are ${r-1\choose 2}+2(r-1)$ such $\{a,b\}$), let $v_{ab}$ be a vertex in a part different from $a,b$ that is distinct from all other vertices being considered.  The number of choices for the $({r+1 \choose 2}+r-2)$-tuple $(v_2, \ldots, v_r, \{v_{ab}\}_{a,b})$ is at least $(1-\od1)c(n, L_{r+1})$. Moreover, the ${r+1 \choose 2}$ edges $e$ and $\{abv_{ab}\}_{a,b}$ form a potential copy of $L_{r+1}$ with $x,y,v_2,\ldots, v_r$ forming the original $K_{r+1}$ whose edges have been expanded.
At least $(\e/2)c(n, L_{r+1})$ of these potential copies of $L_{r+1}$ have a triple from $M$,  otherwise $e$ would lie in at least
$(1-o_{\delta}(1)-\e/2)c(n, L_{r+1})>(1-\e)c(n, L_{r+1})$ copies of $L_{r+1}$ and we
are done.  Suppose that at least $(\e/4)c(n, L_{r+1})$ of these potential copies of $L_{r+1}$ have the triple from $M$ omitting $e$. Since each such triple from $M$ is counted at most $n^{l_{r+1}-6}$ times,  we obtain the contradiction $(\e/4)c(n, L_{r+1})/n^{l_{r+1}-6}\le |M| <\odel(n^3).$
 So at least $(\e/4)c(n, L_{r+1})$ of these potential copies of $L_{r+1}$ have a triple from $M$ containing $x$ or $y$.
Each such triple from $M$ is counted at most $n^{l_{r+1}-5}$ times, so there are at least $(\e/4)c(n, L_{r+1})/n^{l_{r+1}-5}=\e'n^2$ triples from $M$ containing $x$ or $y$ (for suitable $\e'>0$ depending only on $r$).
We may assume by symmetry that $d_M(x) >(\e'/2)n^2$.

We have shown above that for each $e \in B$, there is a vertex $x \in e$ that lies in the (unique) part that has at least two points from $e$, with $d_M(x) >(\e'/2)n^2$. Form a bipartite graph with parts $B$ and $M$, where each $e \in B$ is adjacent to those $f$ in $M$ for which $e \cap f=\{x\}$ and $x$ lies in the part that has at least two points of $e$. Then each vertex of $B$ has degree at least $(\e'/2)n^2$.
Since $|B|>|M|$ we conclude that there exists $f \in M$ adjacent to at least $(\e'/2)n^2$ different $e \in B$
in the way specified above. At least $(\e'/6)n^2$ of these $e \in B$ contain the same point $x \in f$.  Assume wlog that $x \in V_1$.

For each $i \in [r]$ and $\e_1=\e'/100$, define
$$A_i=\{y \in V_i: d_{\cH}(xy)\ge \e_1 n\}.$$

\noindent{\bf Claim.} $|A_i|<\e_1 n$ for some $i \in [r]$.

{\bf Proof of Claim.} Suppose to the contrary that $|A_i|\ge \e_1 n$ for each $i$. Then the number of choices
$(v_1, \ldots, v_r) \in A_1 \times \cdots \times A_r$ is at least $(\e_1 n)^r$. For every pair of distinct vertices $\{a,b\} \subset \{v_1, \ldots, v_r\}$, let $w_{ab}\ne x$ be a vertex in a part different from $a,b$ (there are at least $(1-2/r)n>n/2r$ choices for $w_{ab}$). For every vertex $c \subset \{v_1,\ldots, v_r\}$, let $w_{c}$ be a vertex such that $xcw_c \in \cH$. By definition of $A_i$,  we know that the number of such $w_c$ is at least $\e_1n$. Consequently,
the number of choices for the $(l_{r+1}-1)$-tuple of distinct vertices $(v_1, \ldots, v_r, \{w_{ab}\}_{a,b}, \{w_c\}_c)$ is at least
$$(\e_1 n)^r(n/2r)^{{r \choose 2}}(\e_1n)^r>(\e_1/r)^{r^2}n^{l_{r+1}-1}=\e_2n^{l_{r+1}-1}.$$ Moreover, the ${r+1 \choose 2}$ triples $xcw_c, abw_{ab}$ over all choices of $a,b,c$
 form a potential copy of $L_{r+1}$ with $x, v_1, \ldots, v_r$ forming the original $K_{r+1}$ whose edges have been expanded.
At least $(\e_2/2)n^{l_{r+1}-1}$ of these potential copies of $L_{r+1}$ have a triple from $M$,  otherwise $\#L_{r+1} \ge
(\e_2/2)n^{l_{r+1}-1} >\delta n^{l_{r+1}-1}$
 and we
are done.  Each such triple from $M$ omits $x$ and is therefore counted at most $n^{l_{r+1}-4}$ times (since $x$ is fixed and $L_{r+1}$ has $l_{r+1}$ vertices) so we obtain the contradiction $(\e_2/2)n^{l_{r+1}-1}/n^{l_{r+1}-4}\le |M| <\odel(n^3)$. This completes the proof of the Claim.

Let $B(x)$ be the set of edges of $B$ containing $x$ with at least two vertices in $V_1$.
Then we had earlier shown that $|B(x)|\ge (\e'/6)n^2>10\e_1 n^2$.

 Let $H(x)$ be the set of pairs $\{y,z\}$ such that $xyz \in B(x)$, so one of $y,z \in V_1$ and $|H(x)|=|B(x)|$.
Now $|A_1|\ge \e_1 n$ for otherwise we obtain the contradiction
$$|B(x)| \le \sum_{v \in V_1} d_{H(x)}(v) =\sum_{v \in A_1} d_{H(x)}(v)+\sum_{v \in V_1-A_1} d_{H(x)}(v)\le (\e_1n)n+(n/2)(\e_1n)<2\e_1n^2.$$
The Claim implies that one of $|A_2|, \ldots, |A_r|$ is less than $\e_1n$. By
 symmetry, we may assume that $|A_r|<\e_1n$. The number of edges in $\cH$ containing $x$ and some vertex of $V_r$ is at most $|A_r|n+|V_r|\e_1n<2\e_1n^2$. Hence the number of edges in $B(x)$ that have no vertex in $V_r$ is at least $|B(x)|-2\e_1n^2>8\e_1n^2$.

 Now let us contemplate  moving $x$ from $V_1$ to $V_r$. The edges of $\cH$ containing $x$ whose contribution to $\sum_i ih_i$ decreases (by at most one) must have  a vertex in $V_r$, and their number is at most $2\e_1n^2$. The edges in $B(x)$ that have no vertex in $V_r$ give an increased contribution to $\sum_i ih_i$ (each edge contributes an increase of exactly one), and their number is at least $8\e_1n^2$. All other edges containing $x$
 (i.e. those with $r-1$ vertices in $V_2 \cup \cdots \cup V_{r-1}$) do not change their contribution  to $\sum_i ih_i$.  The net contribution to $\sum_i ih_i$ therefore increases by at least $6\e_1 n^2>0$, thus contradicting the choice of the partition and completing the proof. \qed

\bigskip

\subsection{Exact Counting}

In this subsection we will use Theorem \ref{L1} to prove Theorem \ref{exactLq}.

 \bigskip

\noindent{\bf Proof of Theorem \ref{exactLq}.} Given $q\ge 1$, let $0<\e \le 1/(q+1)$. Then $2(1-\e)c(n, L_{r+1})\ge c(n, L_{r+1})$.   Let $\delta$ and $n_0$ be the outputs of Theorem \ref{Lq} with input $\e$. Choose $n>n_0$ such that it also satisfies $\delta n^{l_{r+1}-1}>q\times c(n, L_{r+1})$ (this is a triviality since $c(n, L_{r+1})=O(n^{l_{r+1}-3})$).

Suppose that $\cH$ is an $n$ vertex 3-graph with $t^3_r(n)+q$ edges. Write $\#L_{r+1}$ for the number of copies of $L_{r+1}$ in $\cH$. Let us prove by induction on $q$ that $\#L_{r+1}\ge q\times c(n, L_{r+1})$. If $q=1$, then  Theorem \ref{L1} and
the definitions of $\e, \delta, n$ imply that
$$\#L_{r+1} \ge  \min\{\delta n^{l_{r+1}-1},\, c(n, L_{r+1}), \,2(1-\e)c(n, L_{r+1})\} \ge c(n, L_{r+1}).$$  Let us assume that $q >1$ and the result holds for  $q-1$.

Let $e_1$ be an  edge of $\cH$ that lies in the maximum number of copies of $L_{r+1}$, say that it lies in $c_1(n)$ copies. If $c_1(n)\ge c(n, L_{r+1})$, then let $\cH_1=\cH-e_1$.  By induction, $\cH_1$ has at least $(q-1)c(n, L_{r+1})$ copies of $L_{r+1}$. These copies are distinct from those containing $e_1$ so we obtain
$$\#L_{r+1} \ge c_1(n)+(q-1)c(n, L_{r+1}) \ge qc(n, L_{r+1})$$ and we are done.

We may therefore assume that $c_1(n)<c(n, L_{r+1})$. Let $e_2$ be an edge of $\cH_1$
that lies in the maximum number $c_2(n)$ of copies of $L_{r+1}$ in $\cH_1$. Since $\cH_1 \subset \cH$, clearly $c_2(n) \le c_1(n)$. Let $\cH_2=\cH_1-e_2$ and continue this process to obtain $e_1, \ldots, e_{q-1}$.  For each $i \le q-1$, Theorem \ref{L1} implies that
$\#L_{r+1}\ge\delta n^{l_{r+1}-1}>qc(n, L_{r+1})$ or
$c_i(n)\ge (1-\e)c(n, L_{r+1})$. In the former case we are done, so we may assume that
$$(1-\e)c(n, L_{r+1}) \le c_{q-1}(n) \le \cdots \le c_1(n) <c(n, L_{r+1}).$$
Consider $\cH_{q-1}=\cH-e_1-e_2 \ldots -e_{q-1}$. Then
$$|\cH_{q-1}|=|\cH|-(q-1)=t^3_r(n)+q-(q-1)=t^3_r(n)+1.$$
Since $c_{q-1}(n)<c(n, L_{r+1})$, Theorem \ref{L1} implies that $\cH_{q-1}$ has at least $2(1-\e)c(n, L_{r+1})$ copies of $L_{r+1}$.  Altogether we have
$$\#L_{r+1} \ge 2(1-\e)c(n, L_{r+1})+\sum_{i=1}^{q-1}c_i(n)\ge (1-\e)(q+1)c(n, L_{r+1})\ge qc(n, L_{r+1})$$
where the last equality follows from $\e\le 1/(q+1)$. This completes the proof. \qed

\section{Concluding Remarks}

$\bullet$ We have given counting results for every triple system for which a stability result is known except for one family which is derived from the expanded cliques.
This was studied in \cite{MP}, and included the triple system $\{123, 145, 167, 357\}$ which is the smallest non-3-partite linear (every two edges share at most one vertex) 3-graph. It appears that our approach will give appropriate counting results for this problem
as well and we did not feel motivated to carry out the details.

 $\bullet$ Our results suggest that whenever one can obtain stability and exact results for an extremal problem, one can also obtain counting results.  However, in each case the argument is different. It would be interesting to unify this approach (at least for certain classes) so one does not have to use new methods for each $F$.  We formulate this as a conjecture.
Say that a 3-graph $F$ is stable if ex$(n,F)$ is achieved uniquely by the $n$ vertex 3-graph $\cH(n)$ for sufficiently large $n$, and every $n$ vertex 3-graph with $(1-o(1))$ex$(n,F)$ edges and no copy of $F$ can be obtained from $\cH(n)$ by changing at most $o(n^3)$ edges.

\begin{conjecture} \label{conj}
Let $F$ be a non 3-partite stable 3-graph. For every positive integer $q$, the following holds for sufficiently large $n$: Every $n$ vertex 3-graph with ex$(n,F)+q$ edges contains at least $qc(n,F)$ copies of $F$.\end{conjecture}

$\bullet$  We have not been able to prove exact counting results for $F_5$ and $B_5$.
The reason for this is that we need to use the minimum degree condition in the proof and we don't know how to get around this technical difficulty.

$\bullet$ All our theorems find $\alpha (1-o(1))n^{\beta}$  copies of $F$ on an edge, or $\delta n^{\gamma}$ copies of $F$ altogether, for suitable $\alpha, \beta, \gamma, \delta$.  However, in each case our proofs  give $\delta n^{\gamma}$ copies of $F$ on a single vertex.

$\bullet$  Our results for $F_5$ appear to be weaker than the other results. In particular, we only allow $q<\delta n$ unlike in the other cases where we allow $q<\delta n^2$.  However, this cannot be improved further.  Indeed, for any $\e>0$ (take $\e=1/2$ for example)
and all $n$, there exists an $n$ vertex 3-graph $\cH$ with $t^3(n)+\e n$ edges and the following two properties:

(1) for every edge $e \in\cH$, the number of copies of $F_5$ containing $e$ is less than $(3-\e)(n/3)^2$

(2) the number of copies of $F_5$ in $\cH$ is less then $\e n^3$.

To see this, let  $T^3(n)$ have parts $X, Y, Z$ and construct $\cH$ as follows. Pick $(x,y) \in X \times Y$, delete $\e n/3$ edges of the form $xyz$ with $z \in Z$, and add $4\e n/3$ edges of the form $x_ixy$ with $x_i \in X$.  Then $|\cH|=t^3(n)+\e n$.  A copy of $F_5$ in $\cH$ must contain an edge $e_i=x_ixy$, and the number of copies containing
$e_i$ is at most $(3-\e)(n/3)^2$. Therefore the total number of copies of $F_5$ in $\cH$ is at most
$(4\e n/3)(3-\e)(n/3)^2<\e n^3$.

$\bullet$ Our results for $L_r$ can be extended to the $k$-uniform case without too much difficulty.  We describe some of the details below.
For $r>k\ge 2$,
 Let $L^k_r$ be the $k$-graph
obtained from the complete graph $K_{r}$ by enlarging each edge
with a set of $k-2$ new vertices.  These sets of new vertices are disjoint for each edge, so $L^k_r$ has $r+(k-2){r \choose 2}$ vertices and ${r \choose 2}$ edges. Write $T^k_r(n)$ for the complete $r$-partite $k$-graph with the maximum number of edges. So $T^k_r(n)$ has vertex partition $V_1 \cup \cdots  \cup V_r$, where $n_i:=|V_i|=\lfloor (n+i-1)/r\rfloor$, and all $k$-sets with at most one point in each $V_i$.  Define
$$t^k_r(n):= |T^k_r(n)| =\sum_{S \in {[r] \choose k}} \prod_{i \in S} n_i.$$
Every set of $r+1$ vertices in $T^k_r(n)$ contains two vertices in the same part, and these two vertices lie in no edge. Consequently, $L^k_{r+1} \not\subset T^k_r(n)$.
The author \cite{M} conjectured, and Pikhurko \cite{P} proved,
that among all $n$ vertex $k$-graphs containing no copy of  $L^k_{r+1}$ ($r\ge k \ge 2$ fixed, $n$ sufficiently large), the unique one with the maximum number of edges is $T^k_r(n)$.  Define $c^k_{r+1}(n)$ to be the minimum number of copies of $L^k_{r+1}$ in a $k$-graph obtained from $T^k_r(n)$ by adding one edge.  The following theorem can be proved by extending the ideas of \cite{P} and Theorem \ref{Lq}'s proof  in the obvious way.

\begin{theorem} \label {Lkq}
Fix $r\ge k \ge 3$. For every $\e>0$ there exists $\delta>0$ and $n_0$ such that the
following holds for $n>n_0$. Let $\cH$ be a $k$-graph with $t^k_r(n)+q$
edges where $q<\delta n^{k-1}$.  Then the number of copies of $L^k_{r+1}$ in
$\cH$ is at least $q(1-\e)c^k_{r+1}(n)$. The expression $q$ is sharp for
$1\le q<\delta n^{k-1}$. Moreover, if the number of copies is less than
$\delta n^{r+(k-2){r+1 \choose 2}}$, then there is a collection of $q$ distinct edges that
each lie in $(1-\e)c^k_{r+1}(n)$ copies of $L^k_{r+1}$ with no two of these edges accounting for the same copy of $L^k_{r+1}$.
\end{theorem}

The exact result for this situation can also be proved using the same methods.

Alon and Pikhurko \cite{AP} proved that ex$(n, L^k(G))=t^k_r(n)$ (for $n>n_0$) where $L^k(G)$ is the $k$-graph  obtained from an $r$-color critical graph $G$
by expanding each edge of $G$ by a new set of $k-2$ vertices. In \cite{KM} we had proved the corresponding counting result for $L^2(G)$ and those ideas combined with the ones in this paper can be used to give similar results for $L^k(G)$.


\end{document}